\documentclass[12pt,leqno]%,dvipdfmx]
{amsart}

\usepackage[margin=1.2in]{geometry}

\usepackage[colorlinks, linkcolor=blue,citecolor=blue, urlcolor=blue]{hyperref}

\usepackage{amsmath,epsfig,graphicx,color, float}
\usepackage{subcaption}
\usepackage{relsize}
\usepackage{comment}

\usepackage{amsmath}
\usepackage{amssymb}
\usepackage{natbib}

\usepackage{enumitem}

\usepackage{algorithm}
\usepackage{algorithmicx}
\usepackage{algpseudocode}

\numberwithin{table}{section}
\numberwithin{figure}{section}
\numberwithin{equation}{section}

\definecolor{darkblue}{rgb}{.2, 0.2,.8}
\definecolor{darkgreen}{rgb}{0,0.5,0.3}
\definecolor{darkred}{rgb}{.8, .1,.1}

\newcommand{\bfY}{\vect{Y}}
\newcommand{\bfy}{\vect{y}}

\newcommand{\bfX}{\vect{X}}

\newcommand{\bfz}{\vect{z}}
\newcommand{\bfZ}{\vect{Z}}

\newcommand{\bfalp}{\vect{\alpha}}

\newcommand{\bfpi}{\vect{\pi}}
\newcommand{\bfT}{\mat{T}}
\newcommand{\bfS}{\mat{S}}
\newcommand{\bfs}{\vect{s}}
\newcommand{\bft}{\vect{t}}
\newcommand{\bfe}{\vect{e}}

\newcommand{\bfI}{\mat{I}}

\newcommand{\bfdelta}{\vect{\delta}}
\newcommand{\bfbeta}{\vect{\beta}}
\newcommand{\bfeta}{\vect{\eta}}

\newcommand{\1}{{\mathbf 1}}
\newcommand{\0}{\mat{0}}

\newcommand{\N}{\mathbb{N}}
\newcommand{\R}{\mathbb{R}}
\newcommand{\E}{\mathbb{E}}
\renewcommand{\P }{{\mathbb P}}

\newcommand{\tod}{\stackrel{d}{\to}}

\newtheorem{lemma}{Lemma}[section]

\newtheorem{theorem}[lemma]{Theorem}

\newtheorem{corollary}[lemma]{Corollary}
\newtheorem{example}[lemma]{Example}

\newtheorem{remark}{Remark}[section]

\usepackage{bm}
\newcommand{\vect}[1]{\pmb{#1}}
\newcommand{\mat}[1]{\boldsymbol{\bm #1}}

\DeclareMathOperator*{\argmax}{arg\,max}

\allowdisplaybreaks

\begin{document}
%\today
%\title[Fitting phase-type frailty models]{ Fitting phase-type frailty models}
\title[Phase-type frailty models]{Phase-type frailty models: A flexible approach to modeling unobserved heterogeneity in survival analysis}
\thanks{The author would like to acknowledge financial support from the Swiss National Science Foundation Project IZHRZ0\_180549.}
%The initial version of this manuscript was developed during my time as a postdoctoral researcher at the University of Bern, where I was financially supported by the Swiss National Science Foundation Project IZHRZ0\_180549.}

\author[J.~Yslas]{Jorge Yslas}
\address{Institute for Financial and Actuarial Mathematics,
University of Liverpool,
L69 7ZL,
Liverpool,
United Kingdom}
\email{Jorge.Yslas-Altamirano@liverpool.ac.uk}

\begin{abstract}
Frailty models are essential tools in survival analysis for addressing unobserved heterogeneity and random effects in the data. These models incorporate a random effect, the frailty, which is assumed to impact the hazard rate multiplicatively. In this paper, we introduce a novel class of frailty models in both univariate and multivariate settings, using phase-type distributions as the underlying frailty specification. We investigate the properties of these phase-type frailty models and develop expectation-maximization algorithms for their maximum-likelihood estimation. In particular, we show that the resulting model shares similarities with the Gamma frailty model, has closed-form expressions for its functionals, and can approximate any other frailty model. Through a series of simulated and real-life numerical examples, we demonstrate the effectiveness and versatility of the proposed models in addressing unobserved heterogeneity in survival analysis.
%Frailty models are survival analysis models which account for heterogeneity and random effects in the data. In these models, the random effect (the frailty) is assumed to have a multiplicative effect on the hazard.
%In this paper, we present frailty models using phase-type distributions as the frailties.
%We explore the properties of the proposed frailty models and derive expectation-maximization algorithms for maximum-likelihood estimation.
%We provide several numerical examples, both for simulated and real-life insurance data.

\end{abstract}
\keywords{frailty model; heavy tails; parameter estimation; phase-type distributions}
\subjclass{Primary 62N02; Secondary 60E05; 60J22}
\maketitle

\section{Introduction}

A core concept in survival analysis is that of the hazard function. This function specifies the instantaneous risk of the event of interest for an individual, given that the individual has not experienced the event previously. Models based on the hazard function are fundamental tools in survival analysis, being Cox's proportional hazards model \citep{cox1972regression} one of the most influential models. This model assumes that the 
ratio of the hazards between any two individuals is constant over time. 
Implicitly, this means that if the model is perfectly specified so that all possible relevant covariates are accounted for, then all individuals in a group with the same covariates have the same risk of the event of interest. 
However, in practice, it is impossible to include all relevant risk factors. 
This unaccounted part is usually known as the unobserved heterogeneity.
The frailty model addresses the problem of unobserved heterogeneity in events. 
The first account of univariate frailty models can be traced back to the work of \cite{beard1959note}, although the term frailty itself was introduced by \cite{vaupel1979impact}. 
This model assumes that an individual's hazard function depends on an unobservable, time-independent random variable known as the frailty.
%A popular choice for the frailty is the Gamma distribution, mainly due to its mathematical tractability and computational feasibility. Another property for the motivation for the use of this specification is that the distribution of the heterogeneity among survivors converges to a Gamma distribution \citep{abbring2007unobserved}. 
%However, no other biological reason makes the Gamma distribution more appealing than other distributions.
%{\red Other popular choices for the law of the frailty, include...}
A popular choice for describing the frailty is the Gamma distribution, a preference driven largely by its mathematical tractability and computational feasibility. One additional rationale for the use of this specification is that the distribution of the heterogeneity among survivors converges to a Gamma distribution \citep{abbring2007unobserved}. Parallelly, another frailty demonstrating mathematical tractability is the positive stable distribution proposed by \cite{hougaard1986class}. The appeal of both the Gamma and positive stable distributions can largely be attributed to the simplicity of their Laplace transforms, which leads to mathematically tractable models. However, it is pertinent to note that this convenience does not translate to satisfactory fitting when applied to real-world datasets. Furthermore, no other biological reasons make these two distributions more appealing than other specifications.
To overcome this limitation, several authors have studied frailty models with diverse frailty distributions exhibiting explicit Laplace transforms. For example, \cite{hougaard1986survival} considered the more general power variance function distribution, which has the Gamma distribution as a particular case and the positive stable distribution as a limiting distribution. 
%However, there is no explicit representation of the density function for this last more general specification.

A phase-type distribution is defined as the distribution of the time until absorption of an otherwise transient time-homogeneous pure-jump Markov process.
 These distributions have a long history going back to the early works of \cite{erlang1909sandsynlighedsregning} and \cite{jensen1954distribution}, and more recently \cite{neuts1975probability}.
 Examples of phase-type distributions include several classical distributions
such as the exponential, Erlang, Coxian, hyperexponential, and any finite mixture of them.
Phase-type distributions have been employed in various domains since they often provide exact, or even explicit, solutions to important problems in complex stochastic models. This is the case, for example, in fields as diverse as biology, operational research, queueing theory, renewal theory, and risk theory, cf., e.g., \cite{bladt2005review,buchholz2014input,hobolth2019phase}. 
Furthermore, the class of phase-type distributions is known to be dense in the set of distributions with support on the positive real numbers in the sense of weak convergence \cite[see][Section 3.2.1]{Bladt2017}. This means that any distribution with support on the positive half-line can be approximated arbitrarily well by a phase-type distribution (of sufficiently large dimension). Moreover, statistical inference for phase-type distributions is a topic well-developed in the literature. For instance, maximum-likelihood estimation was proposed by \cite{asmussen1996fitting} using an expectation-maximization (EM) algorithm, which was subsequently extended in \cite{olsson1996estimation} for the case of censored observations.
% {\red These properties make phase-type distribution attractive modeling alternatives to other flexible tools such as mixture of Gaussians \citep[cf., e.g.,][]{mclachlan2019finite} }
Regarding the use of phase-type distributions in survival analysis, one of the first explorations can be found in \cite{aalen1995phase}, where the author showed the flexibility of the hazard function for phase-type distributions and pointed out that the interpretability of phase-type distributions in terms of absorption times of a hidden Markov chain raises naturally in many survival analysis applications, where the hazard rate of an individual changes depending on different hidden states. For instance, a disease generally goes through various stages of severity. However, as other flexible models, such as mixtures of Gaussians \citep[cf., e.g.,][]{mclachlan2019finite}, phase-type distributions can also be employed as pure modeling tools, and the underlying states do not necessarily have a physical interpretation. 
Later on, in \cite{mcgrory2009fully,tang2012modeling}, phase-type distributions with Coxian structure were considered as survival regression models.
More recently, in \cite{bladt2020survival}, the authors proposed a generalization of the proportional hazards model based on inhomogeneous phase-type distributions \cite[see][]{albrecher2019inhomogeneous} and showed how maximum-likelihood estimation can be performed via an EM algorithm.

In this paper, we propose the use of phase-type distributions as frailties.  
The proposed model's main advantage is that it exploits phase-type distributions' properties, leading to computationally tractable models with desirable characteristics. More specifically, we show that the resulting model has properties that resemble those of the Gamma frailty model, in addition to closed-form expressions for its different functionals.
Moreover, due to the denseness of the class of phase-type distributions, any other frailty model can be approximated by this model.
For application purposes, estimation of these models is essential.  Thus, we proceed to develop an EM algorithm for maximum-likelihood estimation of this model, including the case of right-censoring data and covariates effect. 
Furthermore, we lay out the specific changes to extend the model to the multivariate settings of shared frailty and, more generally, correlated frailty. It is worth noting that, in the same spirit as phase-type distribution, the proposed phase-type frailty models can also have a physical interpretation: the frailty of an individual changes according to different hidden stages. However, the main objective of this paper is to lay out the mathematical foundation of the proposed model and to showcase its modeling capabilities.
%, and not to focus on the interpretability of the model perse. 

The rest of the paper is organized as follows. In Section~\ref{sec:ph}, we present an overview of phase-type distributions and some important properties for our purposes. In Section~\ref{sec:frailty}, we review the univariate frailty model, show how phase-type distributions can be employed as frailties, and derive some relevant properties. In Section~\ref{sec:estimation}, we derive an EM algorithm for maximum-likelihood estimation of the proposed univariate model. 
In Section~\ref{sec:extensions}, we provide multivariate extensions to the shared and correlated frailty cases with an emphasis on their estimation. 
In Section~\ref{sec:examples}, we present several numerical illustrations of practical significance.
 Finally, in Section~\ref{sec:conclusion}, we summarize our findings.

\section{Phase-type distributions}\label{sec:ph}

Let $ ( J_t )_{t \geq 0}$ denote a Markov jump process on the state-space $\{1, \dots, p, p+1\}$, where states $1,\dots,p$ are transient and state $p+1$ is absorbing. In this way, the intensity matrix of $ ( J_t)_{t \geq 0}$ can be written in the form
\begin{align*}
	\mat{\Lambda}= \left( \begin{array}{cc}
		\bfT &  \bft \\
		\0 & 0
	\end{array} \right)\,, 
\end{align*}
where $\bfT $ is a $p \times p$ sub-intensity matrix, representing the transition rates between transient states, and $\bft$ is a $p$-dimensional column vector, denoting the rates of transition into the absorbing state. Additionally, $\0$ is the $p$-dimensional row vector of zeroes. Since the rows of $\mat{\Lambda} $ sum to zero, we have that $\bft =- \bfT \, \bfe$, where $\bfe $ denotes the $p$-dimensional column vector of ones. Assume that the process starts in a transit state $k$ with probability $\pi_k$, that is, $ \pi_{k} = \P(J_0 = k)$, $k = 1,\dots, p$, and let $\bfpi = (\pi_1 ,\dots,\pi_p )$ be the vector of initial probabilities. We further assume that $\P(J_0 = p + 1) = 0$, so that $\sum_{k=1}^p \pi_k = 1$. Then, we say that the time until absorption,
\begin{align*}
	Z = \inf \{ t \geq  0 \mid J_t = p+1 \} \,,
\end{align*}
has a phase-type distribution with representation $(\bfpi,\bfT )$ and we write $Z \sim \mbox{PH}(\bfpi,\bfT )$.  In such a case, $p$ is said to be the dimension of the phase-type distribution. 
We now present an example of phase-type distributions that is particularly relevant for our purposes.

\begin{example}[Generalized Erlang]\rm \label{ex:GErlang}
Let $V_1, \dots, V_p$ be independent random variables with $V_i \sim \mbox{Exp}(\lambda_i)$, $\lambda_i >0$, that is, $\P(V_i > v) = \exp(-\lambda_i v)$, $v>0$, $i = 1,\dots, p$. Then $Z = V_1 + \cdots + V_p$ is phase-type distributed with representation (of dimension $p$)
\begin{align*}
	\bfpi = (1,0,...,0) \,, \quad \bfT= \left( \begin{array}{cccccc}
		-\lambda_1 &  \lambda_1 & 0 & \cdots  & 0 \\
		0 &  -\lambda_2 & \lambda_2  & \cdots  & 0 \\
		0 &  0 & -\lambda_3 & \cdots & 0 \\
		\vdots &  \vdots & \vdots & \ddots & \vdots \\
		0 &  0 & 0 & \cdots  & -\lambda_p \\
	\end{array} \right)\,. 
\end{align*} 

Indeed, $Z$ can be interpreted as the absorption time of a Markov jump process with $p$ transient states as follows. The process initiates in state $1$, sequentially transitioning to subsequent states while holding an exponentially distributed time in each state. This progression continues until reaching state $p$, where it then transitions to an absorbing state. The distribution of $Z$ is called a {\em generalized Erlang} distribution.
 Note that, in particular, the Erlang distribution, that is, the Gamma distribution with an integer shape parameter, stands as a specific instance of phase-type distributions.
%The process starts in state $1$ and jumps to the next state in the sequence (remaining an exponentially distributed time at each state), up to reaching state $p$, from which it jumps to the absorbing state. The distribution of $Z$ is called a {\em generalized Erlang} distribution.
%Note that $Z$ can be obtained by summing the $Z_i$'s in any order. This means that we can always reorder the states to get an alternative (and equivalent) representation of the distribution of $Z$. In other words, the representation is not unique. 
\hfill $\circ$
\end{example}

%\begin{example}[Hyperexponential]\label{ex:HExp} \rm 
%	Consider $V_1, \dots, V_p$ independent random variables, where $V_i \sim \mbox{exp}(\lambda_i)$, $\lambda_i >0$, with density function $f_i(v) = \lambda_i \exp(- \lambda_i v)$, $v>0$, $i = 1,\dots, p$. Let $\pi_1, \dots, \pi_p$ be probabilities such that $\pi_1 + \dots +\pi_p = 1$, and let $f$ be  given by
%	\begin{align*}
%		f(z) = \sum_{i = 1}^{p} \pi_i f_i(z) \,, \, \quad z >0 \,,
%	\end{align*}
%	which corresponds to the density function of a mixture of exponential distributions (also known as {\em hyperexponential}).
%	Then, $f$ has the phase-type representation
%	\begin{align*}
%	\bfpi = (\pi_1,...,\pi_p)\,, \quad \bfT= \left( \begin{array}{ccccc}
%		-\lambda_1 &  0  & \cdots  & 0 \\
%		0 &  -\lambda_2   & \cdots  & 0 \\
%		\vdots  & \vdots & \ddots  & \vdots \\
%		0  & 0 & \cdots  & -\lambda_p \\
%	\end{array} \right)\,. 
%\end{align*}
%	%Note again that any reordering of the states results in the same distribution. Moreover, take, for instance, any exponential distribution with parameter $\lambda>0$. Then, we can represent such an exponential as a mixture of $p$ exponentials with the same parameter $\lambda$ and any arbitrary selection of weights for any $p \in \N$. This exemplifies the fact that phase-type distributions can have equivalent representations of different dimensions. In Section~\ref{sec:estimation}, we will see that the non-uniqueness of the representation of phase-type distributions has to be considered when dealing with their estimation. 
%	\hfill $\circ$
%\end{example}

The density $f_Z$ and survival function $S_Z$ of $Z \sim  \mbox{PH}(\bfpi , \bfT )$ are given by the closed-form expressions
\begin{eqnarray*}
 f_Z(z) &=& \vect{\pi}\exp \left( \mat{T} z \right)\vect{t} \,, \quad z>0\,, \label{eq:dens-IPH} \\
 S_Z(z) &=&  \vect{\pi}\exp \left(  \mat{T} z \right)\vect{e} \,, \quad z>0\,, \label{eq:cdf-IPH}
\end{eqnarray*}
where the exponential of a matrix $\mat{A}$ is defined by the power series
\begin{align*}
	\exp(\mat{A}) = \sum_{i= 0}^{\infty} \frac{\mat{A}^{i}}{i!} \,.
\end{align*}

%Figure~\ref{fig:denph} presents some density shapes for the phase-type particular instances generalized Erlang and hyperexponential.
%
%\begin{figure}[h]
%\centering
%\includegraphics[width=0.49\textwidth]{dgerlang.png}
%\includegraphics[width=0.49\textwidth]{dhyper.png}
%\caption{Density functions of some phase-type models: generalized Erlang (left), and hyperexponential (right).}
%\label{fig:denph} 
%\end{figure}

Not only are there closed-form expressions for the survival and density functions, but this also applies to the Laplace transform and moments. More specifically, the moments are given by 
\begin{align*}
	\E(Z^n) = n! \, \bfpi (-\bfT)^{-n}\bfe  \,, \quad n \in \N \,,
\end{align*}
and the Laplace transform $\mathcal{L}_Z$ can be expressed as
\begin{align}\label{eq:phlaplace}
	\mathcal{L}_Z(u) = \E\left(\exp(-uZ)\right)= \bfpi (u \mat{I} - \bfT )^{-1} \bft = \bfpi (u (-\bfT)^{-1}  + \mat{I})^{-1} \bfe  \,,
\end{align}
where $\mat{I}$ is the identity matrix of appropriate dimension. This Laplace transform is well-defined for $u$ larger than the largest real eigenvalue of $\bfT$ and, in particular, for $u \geq 0$.
 Furthermore, the derivatives of order $n \in \N $ of $\mathcal{L}_Z$ are given by
\begin{align*}
	\mathcal{L}_{Z}^{(n)}(u) = (-1)^{n}\E\left(Z^{n}\exp({-uZ})\right)= (-1)^{n} n! \,\bfpi (u \mat{I} - \bfT )^{-1-n}\, \bft \,.
\end{align*}

In addition to their mathematical tractability, another property that makes phase-type distributions attractive modeling tools is that they form a class that is dense (in the sense of weak convergence) in the set of all distributions on the positive real-line, as stated in the following result which can be found in, e.g., Section 3.2.1 of \cite{Bladt2017}. 

\begin{theorem}\label{theo:weakph}
Let $V$ be any non-negative random variable. Then there exists a sequence of random variables $(Z_n)_{n\geq 1}$, where $Z_n \sim \mbox{PH}(\bfpi_n, \bfT_n)$, $n \geq 1 $, such that 
\begin{align*}
	Z_n \tod V \,, \quad n\to \infty \,,
\end{align*}
where $\tod $ denotes convergence in distribution or weak convergence.
\end{theorem}

%In particular, this means that phase-type distributions can approximate arbitrarily closely any distribution with support on the positive real-line.

\begin{remark}\label{rem:weak} \rm 
	It might be contested that weak convergence does not inherently guarantee the existence of a $Z_n$ with distribution function $F_{Z_n}$ that can approximate the distribution function $F_Y$ of $Y$ to any specified precision, given that weak convergence operates on a pointwise basis. Nonetheless, in the context of a continuous $F_Y$, Lemma 3.2 in \cite{rao1962relations} asserts that weak convergence does, in fact, imply uniform convergence. This means that for any chosen precision parameter $\epsilon >0$, there exists $Z_n$ such that $\sup_{z>0} |F_{Z_n}(z) - F_Y(z)| < \epsilon$, thus achieving the desired approximation within the stipulated precision.
	%One may argue that weak convergence does not mean that there exists $Z_n$ with distribution function $F_{Z_n}$ that approximates the distribution function $F_Y$ of $Y$ at any given precision since weak convergence is only a pointwise type of convergence. However, if $F_Y$ is continuous, Lemma 3.2 in \cite{rao1962relations} states that weak convergence implies uniform convergence in the sense that for any given precision $\epsilon >0$, there exists $Z_n$ such that $\sup_{z>0} |F_{Z_n}(z) - F_Y(z)| < \epsilon$. 
	%Through this paper, most of the introduced models satisfy the continuity condition of the distribution function, and hence, weak convergence can be thought of as uniform convergence. 
\end{remark}

For an extensive and modern treatment of phase-type distributions, we refer to \cite{Bladt2017}.

\section{ The univariate phase-type frailty model}\label{sec:frailty}
\subsection{ Standard univariate frailty models}
Recall that for a continuous random variable $Y$, the hazard function $\mu_Y$ is given by 
\begin{align*}
	\mu_Y(t) = \frac{f_Y(t)}{S_Y(t)} \,.
\end{align*}

A related function is the cumulative hazard function $M_Y$, which is given by 
\begin{align*}
	M_Y(t) = \int_0^{t} \mu_Y (s) ds = -\log (S_Y(t)) \,.
\end{align*}

In the (univariate) frailty model, one assumes that the hazard function of an individual is proportional to an unobservable non-negative random variable $Z$. More specifically, let $Y$ be a random variable with conditional hazard given $Z$ of the form
\begin{align}\label{eq:frailty}
	\mu(t; Z) = Z \mu(t) \,.
\end{align}

In this context, $\mu$ is known as the baseline hazard function, and the random variable $Z$ as the frailty. Note that the conditional survival function of $Y \mid Z = z$ is given by 
\begin{align*}
	S_{Y|Z}(y | z) = \exp \left( -z \int_0^y \mu(t)dt \right) = \exp \left( -z M(y)\right) \,.
\end{align*}

Thus, the survival function of $Y$ is given by
\begin{align*}
	S_Y(y) = \int_0^\infty S_{Y|Z}(y | z) dF_Z(z) = \int_0^\infty \exp \left( -z M(y)\right) dF_Z(z)  = \mathcal{L}_Z(M(y)) \,.
\end{align*}
%where $\mathcal{L}_Z$ is the Laplace transform of $Z$.

Furthermore, the model can be extended to incorporate predictor variables $\bfX = (X_1,\dots, X_h)$ via 
\begin{align*}
	\mu(t; Z, \bfX) = Z \mu(t) \exp ( \bfX \bfbeta) \,,
\end{align*}
where $\bfbeta$ is an $h$-dimensional column vector of regression parameters. 
Consequently, a frailty model is a generalization of the well-known proportional hazards model. 
For the sake of simplicity, we will restrict our treatment to model \eqref{eq:frailty} sometimes to focus on the main ideas of the frailty model.
Next, we present two classic examples of frailties that belong to the more general family of  positive variance functions (PVF) \cite[see, e.g.,][Section 3.7]{wienke2010frailty}.

\begin{example}[Gamma frailty] \rm
	Assume that $Z \sim \mbox{Gamma}(\alpha, \gamma)$, $\alpha, \gamma >0$. Then, the Laplace transform of $Z$ is given by 
	\begin{align*}
	\mathcal{L}_Z(u) = (1 + u / \gamma)^{-\alpha} \,.
\end{align*}

Thus, the survival function of $Y$ takes the explicit representation
\begin{align*}
	S_Y(y)  = (1 + M(y) / \gamma)^{-\alpha}  \,.
\end{align*} 

To avoid over-parametrization, it is common practice to assume that $\E(Z) = 1$. In the Gamma frailty case, this can be done by setting $\alpha = \gamma = 1 /\sigma^2$ , $\sigma > 0$, so that $\E(Z) = 1$ and $\mbox{var}(Z) = \sigma^2$.
\hfill $\circ$
\end{example}

\begin{example}[Inverse Gaussian frailty] \rm 
Consider inverse Gaussian frailty $Z$  with  parameters $\nu, \lambda > 0$ and density function
\begin{align*}
	f_Z(z) = \frac{\sqrt{\lambda}}{\sqrt{2 \pi z^3}} \exp \left( -\frac{\lambda}{2 \nu^2 z} (z - \nu )^2 \right) \,, \quad z>0 \,.
\end{align*}
Then, the Laplace transform of $Z$ is 
\begin{align*}
	\mathcal{L}_Z(u) = \exp \left( -\frac{\lambda \sqrt{1 + 2 \nu ^2 u / \lambda}}{\nu} + \frac{\lambda}{\nu} \right) \,.
\end{align*}

Take the particular case $\nu = 1$ and $\lambda = 1/\sigma^2$, so that $\E(Z) = 1	$ and $\mbox{var}(Z) = \sigma^2$. In this way, 
\begin{align*}
	S_Y(y) =   \exp \left( \frac{1}{\sigma^2} \left(1 - \sqrt{1 + 2 \sigma ^2  M(y)  }\right) \right) \,.
\end{align*}
\hfill $\circ$
\end{example}

 We refer to, e.g., \cite{wienke2010frailty} for a comprehensive account of frailty models.

\subsection{Phase-type frailty}
 We now introduce our novel model by considering phase-type distributions as frailties. More specifically,
let $Z \sim \mbox{PH}(\bfpi, \bfT)$ with Laplace transform \eqref{eq:phlaplace}. Then, the survival function of $Y$ is given by
\begin{align*}
	S_Y(y) =  \bfpi (M(y) \mat{I} - \bfT )^{-1} \bft \,,
\end{align*}
with corresponding density function
\begin{align}\label{eq:phfrailty}
	f_Y(y) =  \mu(y ) \bfpi (M(y) \mat{I} - \bfT )^{-2} \bft \,.
\end{align}

In particular, this implies that the hazard function of $Y$ is given by 
\begin{align}\label{eq:phhazard}
	\mu_Y(y) = \mu(y ) \frac{\bfpi (M(y) \mat{I} - \bfT )^{-2} \bft }{\bfpi (M(y) \mat{I} - \bfT )^{-1} \bft } \,.
\end{align}

\begin{remark}  \rm 
	 It is noteworthy that the Gamma frailty model with integer shape parameters, that is, $\mbox{Gamma}(k, \gamma)$ with $k \in \N$,  is a specific instance of the phase-type frailty model, as delineated in Example~\ref{ex:GErlang}.
	  %Furthermore, we know that for non-negative random variables, weak convergence implies convergence of the Laplace transform. Then, the denseness of the phase-type class of distributions implies weak convergence of the set of phase-type frailty models to any other frailty model, meaning that any frailty model can be approximated arbitrarily well by a phase-type frailty model. 
\end{remark}

The next Corollary shows that the phase-type frailty model inherits the denseness property of phase-type distributions in the following sense: 
\begin{corollary}
Let $Y$ be a random variable whose distribution is given by a frailty model with frailty a non-negative random variable $Z$ and baseline hazard function $\mu$. Then, there exists a sequence of random variables $(Y_n)_{n\geq 1}$, where $Y_n$ is phase-type frailty distributed with frailty $Z_n \sim \mbox{PH}(\bfpi_n, \bfT_n)$ and the same baseline hazard function $\mu$, $n \geq 1$, such that
\begin{align*}
		Y_n \tod Y \,,\quad   n \to \infty \,.
\end{align*}
%	Let $Z$ be a positive random variable and let $Y$ be a random variable with conditional hazard given by \eqref{eq:frailty}. Then there exist $Z_k \sim \mbox{PH}(\bfpi_k, \bfT_k)$, $k \geq 1$, such that $Y_k$ defined as in \eqref{eq:frailty} with this $Z_k$ satisfies that 
%	\begin{align*}
%		Y_k \tod Y \,,\quad   k \to \infty \,.
%	\end{align*}
\end{corollary}

\begin{proof}
For given $Z$, Theorem~\ref{theo:weakph} states that there exists a sequence $(Z_n)_{n\geq 1}$, where $Z_n \sim \mbox{PH}(\bfpi_n, \bfT_n)$, $n \geq 1$, such that $Z_n \tod Z$. Now, using that 
	for non-negative random variables, weak convergence implies pointwise convergence of the Laplace transform, we have that for all $y>0$
	\begin{align*}
		\lim_{n \to \infty}S_{Y_n}(y) = \lim_{n \to \infty}\mathcal{L}_{Z_n}(M(y)) = \mathcal{L}_{Z}(M(y)) =  S_{Y}(y)   \,.
	\end{align*}
	The result follows. 
\end{proof}

In other words, the result above states that any frailty model can be approximated arbitrarily well by a phase-type frailty model (cf. Remark~\ref{rem:weak}).

%\begin{remark}
%It is well-known that phase-type distributions have identifiability issues.
%Namely, different dimension and parameter configurations may lead to very similar or exactly the same density shapes.
% Thus, while one can assume without loss of generality that the mean of the phase-type distribution is one (the class of phase-type distributions is close under scalar constant multiplication), the identifiability issues would preserve. Hence, we do not make any assumptions on the mean of the phase-type frailty here.
%\end{remark}

We now show a series of properties of this model that resemble those of the Gamma frailty model, for instance, those found in Section 3.3 of \cite{wienke2010frailty}: 
 
% {\red Distributions whose density can be written in terms of a matrix-exponential, including the phase-type distributions, are known as matrix-exponential distributions (see  for further details) and consist exactly of the distributions with rational Laplace transform.}
 
i) An interesting property of the phase-type distribution is that the conditional frailty density, conditional on $Y > t$, belongs to the more general family of matrix-exponential distributions, essentially distributions with rational Laplace transform \citep[see][for further details]{Bladt2017}. Indeed, for $t>0$, we have that
\begin{align*}
	f_Z(z\mid Y>t) & = \frac{\exp(-z M(t)) f_Z(z)}{\P(Y>t)}\\
	& =\frac{\exp(-z M(t)) \bfpi \exp({\bfT z})\bft}{\bfpi (M(t) \mat{I} - \bfT )^{-1} \bft} \\
	& = \frac{\bfpi ( M(t) \bfI - \bfT)^{-1} (-\bfT) }{\bfpi (M(t) \mat{I} - \bfT )^{-1} \bft} \exp({(\bfT - M(t) \bfI)  z}) ( M(t) \bfI - \bfT)\bfe \,,
\end{align*}
which corresponds to a matrix-exponential density with a starting vector
\begin{align*}
	\tilde{\bfpi} = \frac{ \bfpi (M(t) \mat{I} - \bfT )^{-1} (-\bfT)}{ \bfpi (M(t) \mat{I} - \bfT )^{-1} \bft} 
\end{align*}
and generator
\begin{align*}
	\tilde{\bfT} = \bfT - M(t) \mat{I} \,.
\end{align*}

ii) It follows from i) that  
\begin{align}\label{eq:expfrailty}
	\E(Z \mid Y > y) = \frac{\bfpi (M(y) \mat{I} - \bfT )^{-2} \bft }{\bfpi (M(y) \mat{I} - \bfT )^{-1} \bft } \,.
\end{align}
Note that this is a decreasing  function that  reaches its maximum $\E(Z)$ at $y = 0$. Then,  \eqref{eq:phhazard} implies that $\mu_Y(y) \leq \mu(y ) \E(Z)$. \\

iii)  A further consequence of i) is that, given that an individual has survived a time $t >0$, the remaining survival time is frailty distributed with matrix-exponential frailty. This can also be seen from the following straightforward calculation 
\begin{eqnarray*}
	\lefteqn{\P(Y> y + t \mid Y>t )}\\ 
	& & = \frac{\P(Y> y + t) }{\P(Y>t )} \\
	& & = \frac{ \bfpi (M(y +t) \mat{I} - \bfT )^{-1} \bft  }{ \bfpi (M(t) \mat{I} - \bfT )^{-1} \bft }\\
	& & = \frac{ \bfpi (M(t) \mat{I} - \bfT )^{-1} (-\bfT)}{ \bfpi (M(t) \mat{I} - \bfT )^{-1} \bft} \times \\
	& & \quad   ((M(y +t) - M(t)) \mat{I} - (\bfT - M(t) \mat{I})^{-1} \left( - (\bfT - M(t) \mat{I} ) \right) \bfe \,.
\end{eqnarray*}
This corresponds to a frailty model with matrix-exponential frailty with parameters $(\tilde{\bfpi}, \tilde{\bfT})$ as in i)
and hazard function $\tilde{\mu}(y) = \mu (y + t)$. \\

iv) Regarding the tail behavior of the phase-type frailty model, this resembles the one of a Gamma frailty model with integer shape parameter. 
More specifically, using Jordan decomposition, it is easy to see that
\begin{align}\label{eq:asymlaplace}
	\mathcal{L}_Z(u) \sim D u ^{-m} \,, \quad u\to \infty \,,
\end{align}
for some positive constant $D$ and $m \in \N$ depending on $\bfpi$ and $\bfT$.

%{\red 
%\begin{remark}
%	If the focus of modeling is in the tail 
%	
%	
%	I think even the paper of the asyptotics with the Gamma is an pro argument to my model, since my asyptotics will remain close to one of the values up to certain point. 
%\end{remark}
%}

%{\red The following example presents a univariate phase-type frailty model that can also be employed in other contexts different from survival analysis. 
%
% \begin{example}
%%This example shows that the phase-type frailty model leads to parsimonious distributions that can be employed in other contexts different from survival analysis. 
%Consider $\mu(y) = \theta y^{\theta-1}$, $\theta >0$, so that
%\begin{align*}
%	S_Y(y) &=  \bfpi ( y^{\theta} \mat{I} - \bfT )^{-1} \bft  = \bfpi ( y^{\theta}  (-\bfT)^{-1}  + \mat{I})^{-1} \bfe \,.
%\end{align*}
%We call this the {\em matrix-Pareto type III} distribution to distinguish it from the matrix-Pareto models introduced in \cite{albrecher2019inhomogeneous} and \cite{albrecher2021cph}. Regarding the tail behavior,
%note that \eqref{eq:asymlaplace} yields
%\begin{align*}
%	S_Y(y) \sim D y^{-m \theta } \,, \quad y \to \infty \,,
%\end{align*}
%with $m \in \N$, meaning that the distribution is regularly varying with index $m\theta$. Thus, this model can be employed to describe data that exhibit heavy tails. 
%\end{example}}

\section{Parameter estimation}\label{sec:estimation}

The EM algorithm \citep{dempster1977maximum} is an iterative method for maximum-likelihood estimation. This algorithm is particularly suitable for situations best described as incomplete-data problems.
It indirectly solves the problem of maximizing the incomplete-data likelihood by alternating between an expectation (E) step, consisting of computing the conditional expectation of the complete log-likelihood given the observed data, and a maximization (M) step, which requires the maximization of the expected log-likelihood found in the E-step.
Given that the phase-type component is not observed for a replication of the phase-type frailty model, we are in an incomplete-data set-up, and the EM algorithm shall be employed.
For such a purpose, we assume that $M(\,\cdot\,; \bfalp)$ is a parametric function depending on some vector $\bfalp$.  

In many applications, a large proportion of the data is either not entirely observed or censored. Here, we consider only the case of right-censoring since it is the most common scenario in the context of this model's applications in survival analysis. However, the instances of left-censoring and interval-censoring can be treated by similar means.  
Let $Y^*$ be a survival time with density \eqref{eq:phfrailty}, and let $C$ be a censoring time. Thus, an observation consists of $(Y, \Delta)$, where $Y= \min(Y^*,C)$ and $\Delta = \1(Y^* \leq C)$. Consider $(y_1, \delta_1), \dots ,(y_N, \delta_N)$ an iid sample from this model, which is also denoted by $(\bfy, \bfdelta)$. In this case, our complete data consists of observations from $(Y, \Delta)$  and the phase-type frailty $Z \sim \mbox{PH}(\bfpi, \bfT)$. Hence, 
 the complete likelihood $L_c$ is given by  
\begin{align*}
	L_c(\bfalp, \bfpi, \bfT ; (\bfy, \bfdelta) ) = \prod_{n=1}^{N} (z_n \mu(y_n; \bfalp))^{\delta_n} \exp(-z_n M(y_n; \bfalp)) f_Z(z_n ; \bfpi , \bfT) \,,
\end{align*}
where $z_n$, $n = 1, \dots, N$, denotes the (not observed) values of $Z$.
Consequently, disregarding the components not reliant on any parameters, the complete log-likelihood  is given by 
\begin{align*}
	l _c(\bfalp, \bfpi, \bfT ; (\bfy, \bfdelta) )  = \sum _{n=1}^{N} \big[ \delta_n\log(\mu(y_n; \bfalp)) - z_n M(y_n; \bfalp) + \log( f_Z(z_n ; \bfpi, \bfT)) \big] \,.
\end{align*}

\textbf{E-step}

Let us denote by $\bfalp^{(k)}$, $\bfpi^{(k)}$, and $\bfT^{(k)}$ the current parameters after $k$ iterations.  For the $(k+1)$-th iteration,
we first require computing the conditional expectation of the log-likelihood given the observed data. 
More specifically, denoting by $\E^{(k+1)} (\cdot)$ the expectation using the current parameters, we require to compute $\E^{(k+1)} (Z \mid (\bfy, \bfdelta) )$ and $\E^{(k+1)} (\log (f_Z(Z ; \bfpi, \bfT)) \mid (\bfy, \bfdelta))$. 
For that purpose, we first consider one (generic) data point ($N=1$) and let $(y, \delta) = (y_1, \delta_1)$. Then, for the non-censored case ($\delta = 1$), we have that
\begin{align*}
	\E^{(k+1)} (Z \mid (y,1))  &= \E^{(k+1)} (Z \mid Y = y)\\ & = \int_0^\infty z f_{Z|Y}(z|y)dz \\
	& =  \int_0^\infty z \frac{f_{Y|Z}(y|z) f_Z(z)}{f_Y(y)} dz \\
	& = \frac{\mu(y)}{ f_Y(y)  } \int_0^\infty z^2 \exp(-z M(y)) f_Z(z) dz\\
	& = \frac{\mu(y)}{ f_Y(y)  } \mathcal{L}_Z^{(2)}(M(y) ) \\
	& =  \frac{2 \bfpi^{(k)} (M(y; \bfalp^{(k)}) \mat{I} - \bfT^{(k)} )^{-3} \bft^{(k)} }{  \bfpi^{(k)} (M(y; \bfalp^{(k)}) \mat{I} - \bfT^{(k)} )^{-2} \bft^{(k)} } \,.
\end{align*}
Now, for the right-censoring case ($\delta = 0$) we have that $ \E^{(k+1)} (Z \mid  (y,0)) = \E^{(k+1)} (Z \mid  Y>y ) $ is given by \eqref{eq:expfrailty}.

Regarding the logarithmic term, we have that
\begin{eqnarray*}
	\lefteqn{ \E^{(k+1)} (\log (f_Z(Z ; \bfpi, \bfT)) \mid Y = y)} \\
	& &= \int \log (f_Z(z ; \bfpi, \bfT)) f_{Z|Y}(z|y)dz \\
	%& &=  \int \log(f_Z(z ; \bfpi, \bfT)) \frac{f_{Y|Z}(y|z) f_Z(z)}{f_Y(y)} dz \\
	& &=  \int \log(f_Z(z ; \bfpi, \bfT)) \frac{z \exp(- z M(y; \bfalp^{(k)}  ))\bfpi^{(k)} \exp(\bfT^{(k)}z)\bft^{(k)}}{ \bfpi^{(k)} (M(y; \bfalp^{(k)}) \mat{I} - \bfT^{(k)} )^{-2} \bft^{(k)}} dz \,,
\end{eqnarray*}
and
\begin{eqnarray*}
	\lefteqn{ \E^{(k+1)} (\log (f_Z(Z ; \bfpi, \bfT)) \mid Y > y ) }\\
	& &= \int \log(f_Z(z ; \bfpi, \bfT))   \frac{ \exp(-z M(y))  f_Z(z)  }{S_Y(y)} dz \\
	& & = \int \log(f_Z(z ; \bfpi, \bfT)) \frac{ \exp(- z M(y; \bfalp^{(k)}  ))\bfpi^{(k)} \exp(\bfT^{(k)}z)\bft^{(k)}}{ \bfpi^{(k)} (M(y; \bfalp^{(k)}) \mat{I} - \bfT^{(k)} )^{-1} \bft^{(k)}} dz \,.
\end{eqnarray*}

In general, closed-form solutions for these integrals are not available. However, numerical approximations can be employed for the subsequent maximization needed in the M-step. 

Finally, for $N > 1$ we sum over $(y_n, \delta_n)$, $n = 1,\dots,N$, in the formulas above.
\\

\textbf{M-step}

Having found the required expectations, the next step is to maximize the conditional expected log-likelihood with respect to the parameters $\bfalp$  and $(\bfpi, \bfT)$,
by which we will obtain updated parameters $\bfalp^{(k +1)} $  and $(\bfpi^{(k +1)} , \bfT^{(k +1)} )$.  
This  will be done separately.

In full generality, for the parameter $\bfalp$, we write
\begin{align*}
	\bfalp^{(k +1)} = \argmax_{\bfalp} \sum _{n=1}^{N} \left[ \delta_n \log(\mu(y_n; \bfalp)) - \E^{(k+1)}  (Z \mid  (y_n, \delta_n))  M(y_n; \bfalp) \right] \,,
\end{align*}
which can be evaluated numerically.

Regarding the parameters $(\bfpi, \bfT)$  of the phase-type frailty, we will employ the EM algorithm introduced by \cite{asmussen1996fitting} as follows. Recall that such an algorithm can be employed to fit a phase-type distribution to a  given theoretical distribution. More specifically,  let us denote by $f_{(\bfpi, \bfT)}$ the density of the phase-type and $g$ the density of the  given theoretical distribution. By fitting $f_{(\bfpi, \bfT)}$ to $g$, one minimizes the Kullback-Leibler information, or equivalently, maximizes $\int \log(f_{(\bfpi, \bfT)}(w)) g(w) dw $ with respect to $(\bfpi, \bfT)$. Now note that  maximizing 
\begin{align*}
	\E^{(k+1)}  (\log (f_Z(Z ; \bfpi, \bfT)) \mid (\bfy, \bfdelta)) 
	&  =  \sum_{n = 1 }^{N} \Bigg[ \delta_n \int \log(f_Z(z ; \bfpi, \bfT)) f_{Z|Y}(z|y_n) dz 
	\\  &  \qquad  (1 - \delta_n) \int \log(f_Z(z ; \bfpi, \bfT)) \frac{ \exp(-z M(y_n))  f_Z(z)  }{S_Y(y_n)} dz \Bigg] \,
\end{align*}
with respect to $(\bfpi, \bfT)$  is equivalent to maximizing
\begin{eqnarray*}
	\lefteqn{\frac{1}{N}\E^{(k+1)}  (\log (f_Z(Z ; \bfpi, \bfT)) \mid (\bfy, \bfdelta))} \\
	& & =  \int \log(f_Z(z ; \bfpi, \bfT)) \left(\frac{1}{N} \sum_{n=1}^{N} \left[ \delta_n f_{Z|Y}(z|y_n) + (1- \delta_n)\frac{ \exp(-z M(y_n))  f_Z(z)  }{S_Y(y_n)} \right] \right) dz \,.
\end{eqnarray*}
Moreover, 
\begin{align}\label{eq:mixfit}
	\frac{1}{N}\sum_{n=1}^{N} \left[  \delta_n f_{Z|Y}(z|y_n)  + (1- \delta_n)\frac{ \exp(-z M(y_n))  f_Z(z)  }{S_Y(y_n)} \right]
\end{align}
is a density function. 
Thus, this part of the M-step reduces to fitting a phase-type distribution to the density \eqref{eq:mixfit}. In general, the EM algorithm does not necessarily require finding estimators that maximize the complete likelihood to guarantee that the incomplete likelihood increases. It suffices to find parameters such that the complete likelihood increases. 
This is also known as the Generalized EM algorithm \citep[cf.][]{mclachlan2007algorithm}.
Hence, we can apply as many iterations as  necessary of the EM algorithm by \cite{asmussen1996fitting} to obtain our updated parameters $(\bfpi^{(k +1)} , \bfT^{(k +1)} )$. 
 In this way, the EM algorithm described above ensures that 
the likelihood increases at each iteration. 
Therefore, we will have convergence to some
limit, which could be a global maximum but also a local maximum or saddle point \citep{wu1983convergence}. If the latter case occurs, a search for more suitable initial parameters should be performed.
We summarize into a complete and detailed routine, which can be found in Algorithm~\ref{alg:PHfrailty}.

\begin{algorithm}[]
\caption{EM algorithm for the univariate phase-type frailty model}\label{alg:PHfrailty}
\begin{algorithmic}
\normalsize 
\State \textit{\textbf{Input}: Sample $(y_1, \delta_1), \dots ,(y_N, \delta_N)$.}\\
\begin{enumerate}[label=\arabic*.]
\item[0.] Initialize with some ``arbitrary" $( \bfalp^{(0)}, \bfpi^{(0)} , \bfT^{(0)})$ and let $k = 0$.
	\item \textit{E-step:}  Compute 
\begin{align*}
	\E^{(k+1)} (Z \mid (y_n,\delta_n)) 
	& =  \delta_n \frac{2 \bfpi^{(k)} (M(y_n; \bfalp^{(k)}) \mat{I} - \bfT^{(k)} )^{-3} \bft^{(k)} }{  \bfpi^{(k)} (M(y_n; \bfalp^{(k)}) \mat{I} - \bfT^{(k)} )^{-2} \bft^{(k)} } \\
	& \quad  + (1 - \delta_n) \frac{\bfpi^{(k)} (M(y_n; \bfalp^{(k)}) \mat{I} - \bfT^{(k)} )^{-2} \bft^{(k)} }{\bfpi^{(k)} (M(y_n;\bfalp^{(k)}) \mat{I} - \bfT^{(k)} )^{-1} \bft^{(k)} } \,
\end{align*}
for all $n = 1,\dots,N$, and let
\begin{align*}
	 g (z )
	&  = \frac{1}{N} \sum _{n=1}^{N}\Bigg[ \delta_n  \frac{z \exp(- z M(y_n; \bfalp^{(k)}  ))\bfpi^{(k)} \exp(\bfT^{(k)}z)\bft^{(k)}}{ \bfpi^{(k)} (M(y_n; \bfalp^{(k)}) \mat{I} - \bfT^{(k)} )^{-2} \bft^{(k)}} \\
	&  \qquad \qquad\quad + (1 - \delta_n)  \frac{ \exp(- z M(y_n; \bfalp^{(k)}  ))\bfpi^{(k)} \exp(\bfT^{(k)}z)\bft^{(k)}}{ \bfpi^{(k)} (M(y_n; \bfalp^{(k)}) \mat{I} - \bfT^{(k)} )^{-1} \bft^{(k)}} \Bigg] \,.
\end{align*}
	
	\item \textit{M-step: Let} 
	\begin{align*}
	\bfalp^{(k +1)} = \argmax_{\bfalp} \sum _{n=1}^{N} \big[ \delta_n \log(\mu(y_n; \bfalp)) - \E^{(k+1)}  (Z \mid  (y_n, \delta_n))  M(y_n; \bfalp) \big] \,,
\end{align*}
assign $\bfeta:= \bfpi^{(k)}$, $\mat{S}:= \bfT^{(k)}$, $\mat{s}:= \bft^{(k)}$, and iterate as necessary the following routine:

\begin{enumerate}
	\item Calculate
\begin{gather*}
	\hat{\eta}_{i} = \int_{0}^{\infty} \dfrac{\eta_{i} \bfe_{i}^{\top} \exp({\bfS y}) \bfs} { \bfeta \exp({\bfS y}) \bfs} g(y)dy \,,\\[4mm]
	\hat{s}_{il} = \dfrac{ \mathlarger{ \int_{0}^{\infty} \dfrac{s_{il}}{ \bfeta \exp({\bfS y}) \bfs } {\int_{0}^{y}  \bfe_{l}^{\top}   \exp({\bfS (y-u)}) \bfs \bfeta \exp({\bfS u}) \bfe_{i} du} \,g(y)dy}}{\mathlarger{\int_{0}^{\infty} \dfrac{1}{ \bfeta \exp({\bfS y}) \bfs}{\int_{0}^{y} \bfe_{i}^{\top} \exp({ \bfS (y-u)}) \bfs \bfeta \exp({\bfS u}) \bfe_{i} du} \, g(y) dy}} \,,\\[4mm]
	\hat{s}_{i} = \dfrac{ \mathlarger{\int_{0}^{\infty} {s_{i}} \dfrac {\bfeta \exp({\bfS y}) \bfe_{i}} {\bfeta \exp({\bfS y}) \bfs } \, g(y) dy}}{\mathlarger{\int_{0}^{\infty} \dfrac{1}{\bfeta \exp({\bfS y}) \bfs }{\int_{0}^{y} \bfe_{i}^{\top} \exp({\bfS (y-u)}) \bfs \bfeta \exp({\bfS u}) \bfe_{i} du} \, g(y) dy}}\,,\\[4mm]
	\hat{s}_{ii} = - \sum_{l\neq i} \hat{s}_{il} - \hat{s}_i \,,
\end{gather*}
where $\bfe_{i}$ denotes the $i$-th canonical basis vector in $\R^p$.

Let $\hat{\bfeta} = (\hat{\eta}_{1}, \dots, \hat{\eta}_{p})$, $\hat{\bfS} = \{\hat{s}_{il}\}_{i,l = 1,\dots,p}$, and $\hat{\bfs} = (\hat{s}_1, \dots, \hat{s}_p)^{\top}$.
\item Assign $\bfeta:= \hat{\bfeta}$, $\mat{S}:= \hat{\mat{S}}$, $\mat{s}:= \hat{\mat{s}}$  and GOTO (a).
\end{enumerate}\

	\item Assing  $ \bfpi^{(k + 1)}:=\bfeta$, $\bfT^{(k+1)}:= \mat{S} $, $\bft^{(k+1)}:= \mat{s}$, $k := k+1$, and GOTO 1 until a stopping rule is satisfied. 
\end{enumerate}
\State \textit{\textbf{Output}: Fitted parameters $(\bfalp^{(k)}, \bfpi^{(k)} , \bfT^{(k)})$.}
\end{algorithmic}
\end{algorithm}

\begin{remark}[On the structure of the phase-type parameters]\label{rm:structure}\rm
	A notable characteristic of Algorithm~\ref{alg:PHfrailty} is the retention of zero entries in the vector of initial probabilities and the sub-intensity matrix in successive iterations. 
	%This property ensures that if we start with a specific distributional structure, like the generalized Erlang, the resulting fitted distribution retains that form.
	This property is advantageous to obtain a specific structure of the fitted distribution, which, in particular, can help to reduce the number of estimated parameters. For instance, by taking a generalized Erlang structure of the initial parameters, the resulting fitted distribution will retain that form. Below, we present two additional structures that can be employed in practice:
	\begin{itemize}
		\item \textit{Coxian}. This is an extension of the generalized Erlang distribution. In this case, the underlying Markov process can reach the absorbing state from any transient state. Mathematically, the structure is defined by:
		%This structure generalizes the generalized Erlang distribution by allowing the Markov process to reach absorption from any of the transient states. In other words, the representation is as follows
		\begin{align*}
	\bfpi = (1,0,...,0) \,, \quad \bfT= \left( \begin{array}{cccccc}
		-\lambda_1 &  \lambda_1 r_1 & 0 & \cdots  & 0 \\
		0 &  -\lambda_2 & \lambda_2 r_2  & \cdots  & 0 \\
		\vdots &  \vdots & \vdots & \ddots  & \vdots \\
		0 &  0 & 0 & \cdots  & -\lambda_p \\
	\end{array} \right)\,,
\end{align*}
		where $\lambda_i>0$, $r_i \in [0,1]$, $i = 1,\dots, p$.
		\item \textit{Generalized Coxian}.  This structure is similar to the Coxian in terms of the sub-intensity matrix. However, there is added flexibility in terms of the starting state of the Markov process, which is now allowed to be any of the transit states. Specifically, the vector of initial probabilities is:
		%This structure preserves the same form of the sub-intensity matrix of the Coxian structure, but the Markov process is now allowed to start at any state, that is, the vector of initial probabilities is now
		\begin{align*}
	\bfpi = (\pi_1,...,\pi_p) \,,
\end{align*}
with $\pi_i\geq 0$, $i = 1,\dots, p$, $\sum_{i=1}^p \pi_i = 1$.
	\end{itemize}

%	 Figure~\ref{fig:denph2} presents some density shapes for the Coxian and generalized Coxian structures.
%
%\begin{figure}[h]
%\centering
%\includegraphics[width=0.49\textwidth]{dcoxian.png}
%\includegraphics[width=0.49\textwidth]{dgcoxian.png}
%\caption{Density functions of some phase-type models: Coxian (left), and generalized Coxian (right).}
%\label{fig:denph2} 
%\end{figure}

\end{remark}	
 	
%\begin{remark}[On the selection of the phase-type frailty] \rm 
%	The literature in dimension selection for phase-type distributions is rather scarce, and information criteria such as BIC and AIC can over-penalize due to the identifiability issues of phase-type distributions.
%	Hence, a standard approach for selecting the number of phases and the structure of the phase-type is by trial and error. Typically, one starts with low dimensions and sparse structures of the phase-type, compares the results, and assesses the benefit or cost of adding additional dimensions or more general structures by, for example, looking at the changes on the log-likelihood.  
%	%This approach is a similar to other 
%	%much in the same manner as is done for the tuning of hyperparameters of a machine learning method, or the selection of the correct combination of covariates in a linear regression. Such approach is also taken presently.
%		The number of iterations is typically chosen in such a way that changes in subsequent log-likelihoods are negligible. This can be done by considering a number of iterations large enough or by using a stopping rule, such as relative changes in the log-likelihood below a certain given level.
%	In addition to the changes in the log-likelihood, the quality of the fit can also be assessed using visual aids such as QQ-plots.
%	
%	%Tools such as parametric bootstrap are computational costly and the quality of the resulting estimators . Since at each iteration we are not sure which representation is obtained. 
%	
%	
%\end{remark}

%{\red
\begin{remark}[On the computational challenges of the model estimation]\label{rm:comp}\rm
	 Algorithm~\ref{alg:PHfrailty} inherits the computational challenges of the conventional EM algorithm for phase-type distribution introduced in \cite{asmussen1996fitting}. Firstly, like any EM algorithm, converging to a local maximum or saddle point is possible. This issue can be addressed by considering multiple random initializations of the parameters until an adequate model is identified. However, this may require several trials. Secondly, the number of iterations needed for convergence depends also on the number of parameters of the phase-type frailty. Generally, a higher dimensional distribution requires more iterations. Furthermore, as the dimension increases, the computational cost of each iteration also rises, primarily due to the calculation of matrix exponentials and integrals involving matrix exponentials, which are known to be computationally demanding. Nevertheless, this challenge can be mitigated by utilizing sparse structures for the parameters, such as those mentioned in Remark~\ref{rm:structure}, as they typically require fewer iterations to achieve convergence compared to a general phase-type structure of the same order \citep[c.f.][]{asmussen1996fitting}. In addition to these challenges, Algorithm~\ref{alg:PHfrailty} introduces two additional computational costs beyond those encountered in the standard phase-type case. The first cost comes from fitting a phase-type model to the density \eqref{eq:mixfit}, which is typically done by generating a weighted sample from that density. The second cost is due to the additional numerical maximization required to estimate the parameters of the baseline hazard function.
	
	%Given these challenges, we suggest following the commonly applied approach for model selection of phase-type distributions. 
	Given that these challenges closely resemble those encountered in standard phase-type fitting, we recommend following the typical model selection approach used in that context. 
	One should begin by fitting models with phase-type frailties of low dimensions and sparse parameters' structures. Then, gradually consider models with higher dimensions or more general parameters' structures and assess the benefit or cost of this added complexity. This assessment can be done, for example, by analyzing log-likelihood changes to determine whether fit improvements justify the added complexity. Additionally, one can also rely on visual aids such as QQ-plots to further evaluate model adequacy. 

\end{remark}

%}

\begin{remark}\rm 
When considering predictor variables $\bfX$, we have that the conditional density and survival function of $Y | Z = z$ are 
\begin{eqnarray*}
	f_{Y|Z}(y| z; \bfX) & = &  z \exp(\bfX \bfbeta) \mu(y;\bfalp) \exp(- z \exp(\bfX \bfbeta)M(y;\bfalp) )\,,\\
	S_{Y|Z}(y| z; \bfX) & = &  \exp(-z\exp(\bfX \bfbeta)M(y;\bfalp)  ) \,.
\end{eqnarray*}
Then, the complete log-likelihood is given by 
\begin{align*}
	l _c(\bfalp, \bfbeta, \bfpi, \bfT ; (\bfy, \bfdelta), \tilde{\bfX} ) & = \sum _{n=1}^{N} \Big[ \delta_n(\bfX_n \bfbeta + \log(\mu(y_n; \bfalp))) - z_n \exp(\bfX_n \bfbeta) M(y_n; \bfalp) \\
		& \quad \quad \qquad + \log( f_Z(z_n ; \bfpi, \bfT)) \Big] \,,
\end{align*}
where $ \tilde{\bfX} $ denotes the set of observed covariates $\{ \bfX_1, \dots, \bfX_N\}$.
 In this scenario, the only change in the E-step is that the conditional expectations should be modified accordingly to include the factor $\exp(\bfX \bfbeta)$ in the formulas. Regarding the M-step, the maximization of $\bfalp$ and $\bfbeta$ should be done conjunctly. More specifically, we need to compute 
\begin{align*}
	(\bfalp^{(k + 1)}, \bfbeta^{(k + 1)}) & = \argmax_{(\bfalp, \bfbeta)} \sum _{n=1}^{N} \Big[ \delta_n(\bfX_n \bfbeta + \log(\mu(y_n; \bfalp))) \\
	& \quad \quad \quad \quad \quad \qquad - \E^{(k + 1)} (Z \mid (y_n, \delta_n), \bfX_n) \exp(\bfX_n \bfbeta)   M(y_n; \bfalp) \Big] \,.
\end{align*} 
\end{remark}

\section{Multivariate extensions}\label{sec:extensions}

\subsection{Shared phase-type frailty} As a first multivariate extension, we show how phase-type distributions can be employed in the shared frailty model. In the shared frailty model, one assumes that a group of individuals, $\bfY = (Y_1, \dots, Y_d)$, is conditionally independent given the frailty and that their hazard functions are of the form 
\begin{align*}
	\mu_j(t; Z, \bfX_j) = Z \mu_j(t) \exp(\bfX_j \bfbeta) \,, 
\end{align*}
where $\mu_j$ are baseline hazard functions and $\bfX_j$ are $h$-dimensional row vectors of covariates, $j= 1,\dots, d$. In this way, the conditional joint survival function of $\bfY|Z=z$ (omitting covariates) is given by 
\begin{align*}
	S_{\bfY | Z}(\bfy | z) = \prod_{j = 1 }^{d} \exp(- z M_j(y_j)) = \exp \left( -z \sum_{j =1}^{d} M_j(y_j) \right) \,.
\end{align*}
%with corresponding conditional joint density
%\begin{align*}
%	f_{\bfY | Z}( \bfy | z) = \prod_{j = 1 }^{d} z \mu_j(y_j) \exp( - z M_j(y_j) ) \,.
%\end{align*}

Then, the joint survival function of $\bfY$ is given by 
\begin{align*}
	S_{\bfY }(\bfy ) = \mathcal{L}_Z \left( \sum_{j =1}^{d} M_j(y_j) \right) \,,
\end{align*}
with corresponding joint density
\begin{align*}
	f_{\bfY}( \bfy ) = \left( \prod_{j =1}^{d} \mu_j(y_j) \right) \mathcal{L}_Z^{(d)} \left( \sum_{j =1}^{d} M_j(y_j) \right)\,.
\end{align*}

In the phase-type particular case, we have that 
\begin{align*}
	S_{\bfY }(\bfy ) = \bfpi  \left( \sum_{j =1}^{d} M_j(y_j) 
\mat{I} - \bfT \right)^{-1} \bft \,,
\end{align*}
and 
\begin{align*}
	f_{\bfY }(\bfy ) = d!\left( \prod_{j =1}^{d} \mu_j(y_j) \right) \bfpi  \left( \sum_{j =1}^{d} M_j(y_j) 
\mat{I} - \bfT \right)^{-1 - d} \bft \,.
\end{align*}

 We now give some insight into the dependence structure of the proposed model. For such a purpose, we limit ourselves to the bivariate case. 
  A standard way to study the dependence structure in multivariate frailty models is through the cross-ratio function  introduced by \cite{clayton1985multivariate}. This is a measure of local dependence describing changes over time, and it is defined as follows
 \begin{align*}
	\theta(y_1,y_2) &= \frac{\mu(y_1 \mid Y_2 = y_2)}{ \mu(y_1 \mid Y_2 >y_2)}\,.
\end{align*}

 This function can also be written in terms of the joint survival function  and its derivatives as \citep[cf.][]{oakes1989bivariate} 
\begin{align*}
	\theta(y_1,y_2) &= \frac{S_{\bfY }(y_1, y_2) \frac{\partial^2}{\partial y_1 \partial y_2}S_{\bfY }(y_1, y_2)}{ \frac{\partial}{\partial y_1}S_{\bfY }(y_1, y_2) \frac{\partial}{\partial y_2}S_{\bfY }(y_1, y_2)}\,.
\end{align*}

Thus, in the phase-type frailty case, we obtain
\begin{align*}
	 \theta(y_1,y_2)  = \frac{2 \bfpi \left( (M_1(y_1) + M_2(y_2)) \mat{I} - \bfT \right)^{-1} \bft \bfpi \left( (M_1(y_1) + M_2(y_2)) \mat{I} - \bfT \right)^{-3} \bft}{(\bfpi \left( (M_1(y_1) + M_2(y_2)) \mat{I} - \bfT \right)^{-2} \bft ) ^2 } \,.
\end{align*}

\begin{figure}[h]
\centering
\includegraphics[width=0.49\textwidth]{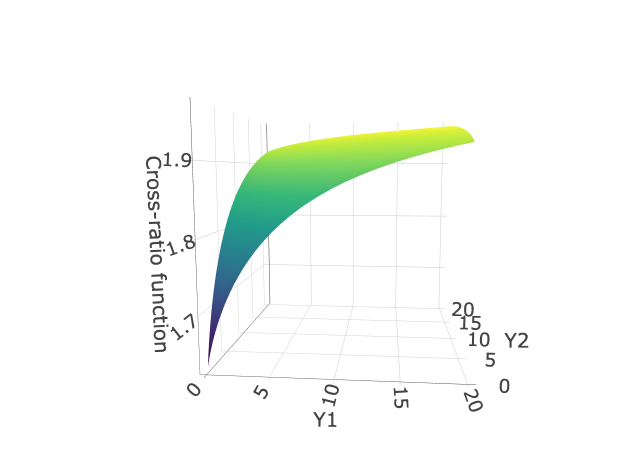}
\includegraphics[width=0.49\textwidth]{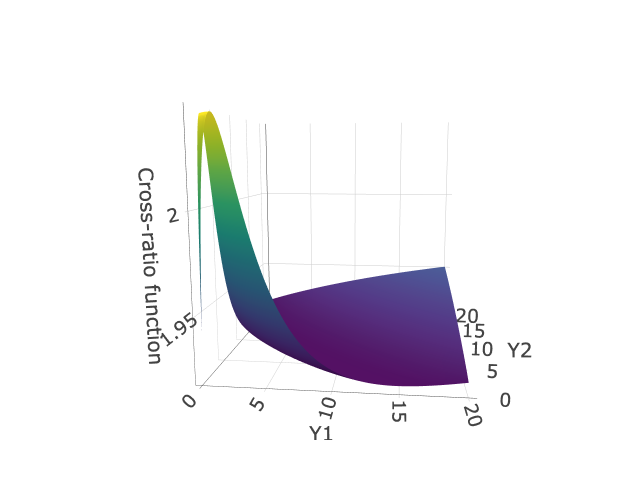}
\caption{Cross-ratio function of a shared phase-type frailty model with frailty a mixture of an exponential and an Erlang distribution (left), and cross-ratio function of a shared phase-type frailty model with frailty a generalized Coxian distribution (right). Weibull baseline hazards were employed in both cases.}
\label{fig:cross}
\end{figure}

\begin{figure}[h]
\centering
\includegraphics[width=0.49\textwidth]{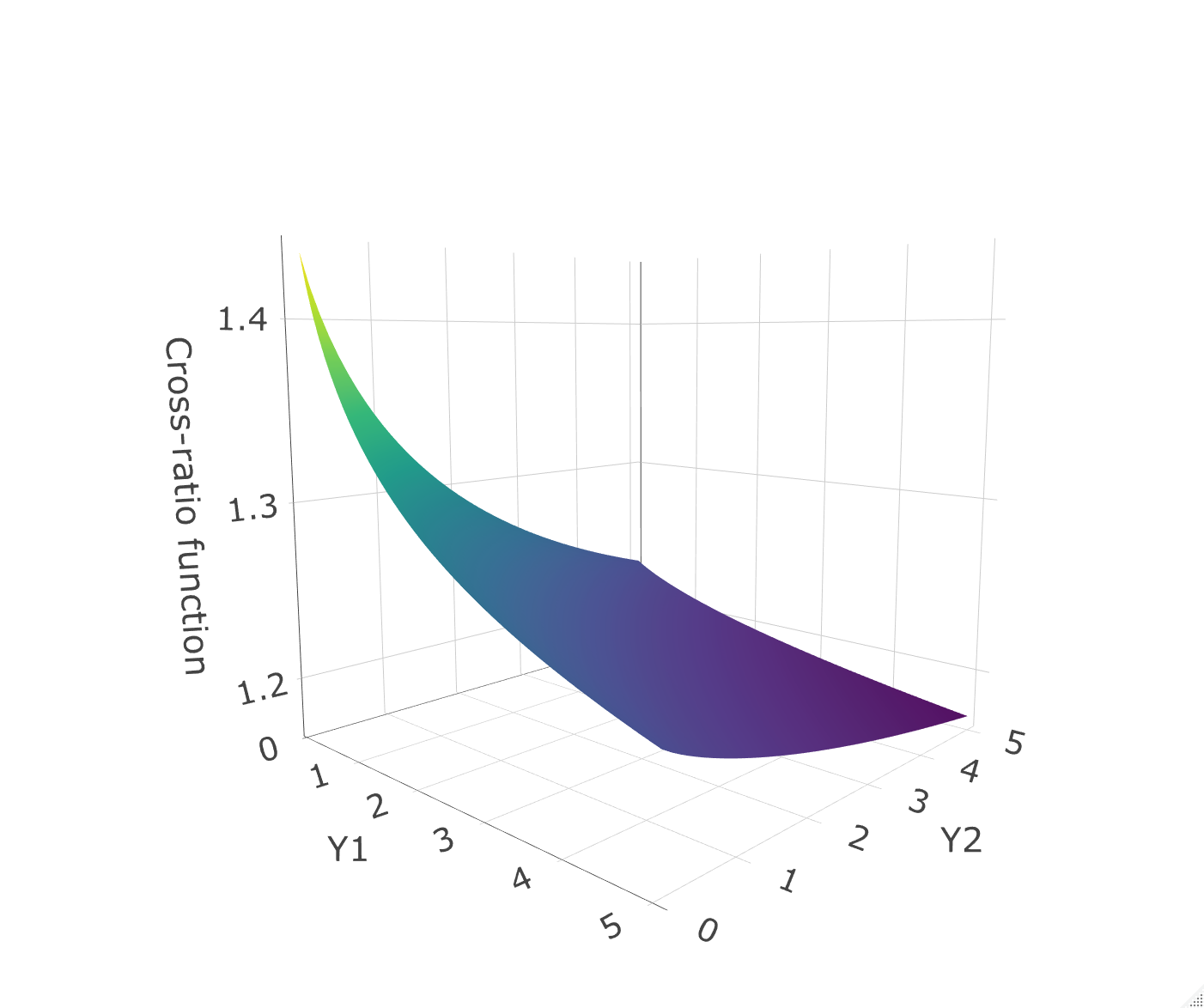}
\includegraphics[width=0.49\textwidth]{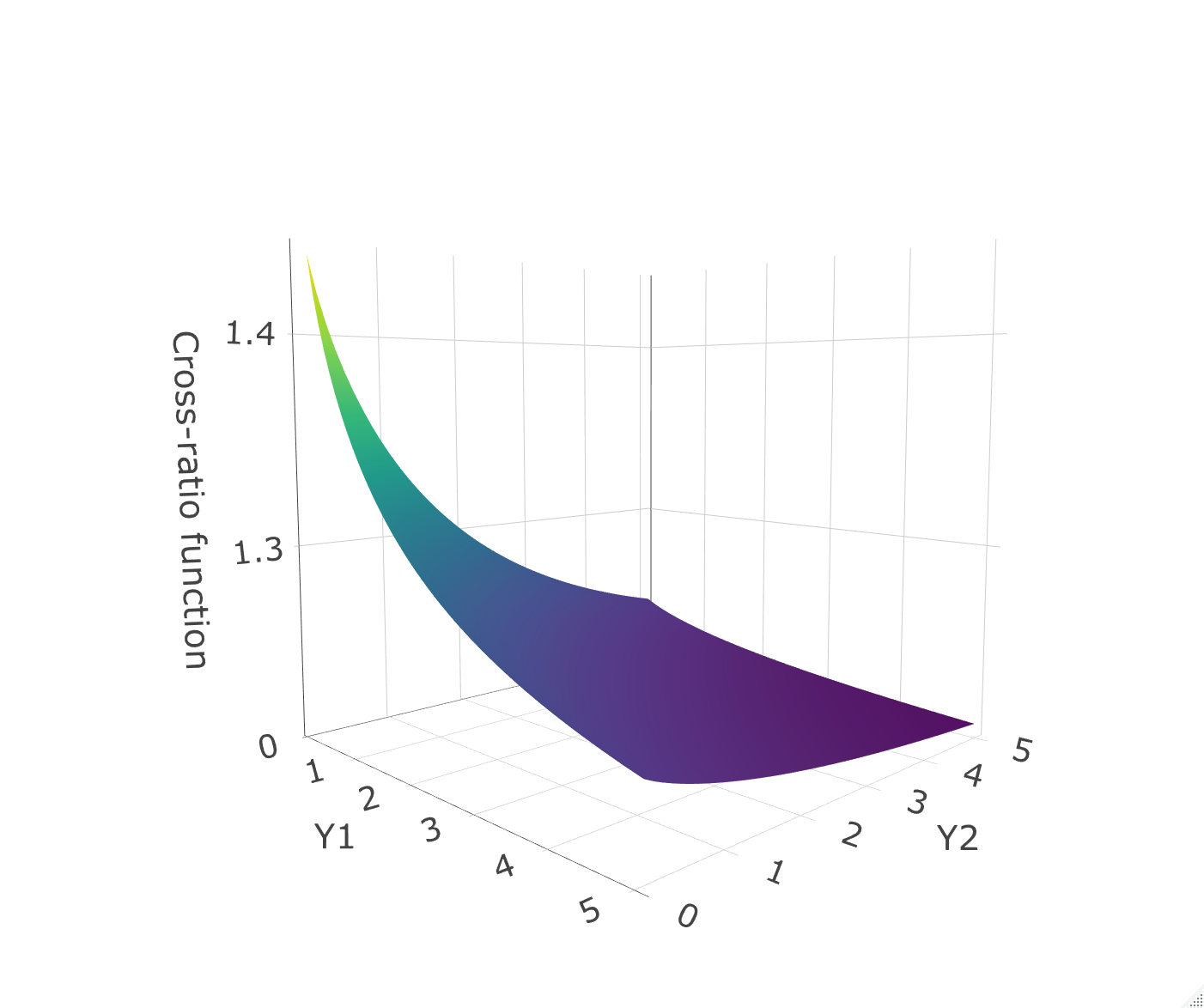}
\caption{Cross-ratio function of a shared inverse Gaussian frailty model (left), and cross-ratio function of a shared phase-type frailty model that approximates the shared inverse Gaussian model (right). Weibull baseline hazards were employed in both cases.}
\label{fig:cross2}
\end{figure}

In stark contrast to the shared Gamma frailty model, where the cross-ratio function is always constant, the above function can have very different behaviors (including also a constant behavior). This is illustrated in Figure~\ref{fig:cross}, where we plot the cross-ratio function of two distinct shared phase-type frailty models. However,
it is easy to see that 
\begin{align*}
	\theta(y_1,y_2)  \to 1 + \frac{1}{m} \,, \quad y_1, y_2 \to \infty\,, %\quad m \in \N \,.
\end{align*}
with $m \in \N$ as in \eqref{eq:asymlaplace}.
Note again that the denseness of the phase-type class of distributions allows us to approximate any other shared frailty model. As an example, Figure~\ref{fig:cross2} shows the cross-ratio function of a shared frailty model with inverse Gaussian frailty, which is then approximated by the corresponding function of a shared phase-type frailty model. 
%Interestingly, Figure~\ref{fig:cross} and Figure~\ref{fig:cross2} illustrate that by using phase-type distribution, the cross-ration function can have very distinct forms, including constant, decreasing, increasing, and "bathtub-shape." This is of particular relevance, as early on, \cite{paik1994multivariate} stated that no frailty models are known to yield bathtub-shaped cross-ratio functions. Later on, \cite{paddy2012relative} showed that bathtub-shapes can be obtained using mixed distributions that have an atom at zero and are continuous on the positive real line. Nevertheless, the authors also highlighted the need for other frailty distributions capable of generating more flexible shapes, particularly without having an atom at zero. 
Interestingly, Figures~\ref{fig:cross} and ~\ref{fig:cross2} demonstrate that with the use of phase-type distribution, the cross-ratio function can adopt a multitude of distinct shapes. This includes patterns that are constant, decreasing, increasing, and even ``bathtub-shaped." This finding holds notable significance, especially when considering the prior observation by \cite{paik1994multivariate}, which posited that no known frailty models produced bathtub-shaped cross-ratio functions. Subsequent work by \cite{paddy2012relative} established that such bathtub-shapes can indeed be derived by using, for example, mixed distributions -- specifically those with an atom at zero and are continuous on the positive real line. However, the same study also underscored the prevailing need for innovative frailty distributions that can craft more flexible cross-ratio functions without necessarily incorporating an atom at zero. In this light, phase-type distributions emerge not just as a novel approach but as a robust solution to this problem.

Another important measure of dependence is the upper tail dependence coefficient $\lambda_U $, defined as 
\begin{align*}
	\lambda_U = \lim_{q \to 1^{-}}\P(Y_1> F^{\leftarrow}_{Y_1}(q) \mid Y_2> F^{\leftarrow}_{Y_2}(q) )	\,,
\end{align*}
where $F^{\leftarrow}_{Y_j}$ denotes the generalized inverse of $F_{Y_j}$, $j = 1,2$. 
For shared frailty models, the above expression takes the particular form
\begin{align*}
	\lambda_U = \lim_{q \to 1^{-}}\frac{\mathcal{L}_Z(2 \mathcal{L}_Z^{-1}(1 - q))}{1 - q}	\,.
\end{align*}

Then, the specific expression of $\mathcal{L}_Z$ for phase-type distributions yields $\lambda_U = 2^{-m}$,  $m \in \N$ .

\subsubsection{Estimation}
We now show how the estimation of the shared phase-type frailty model can be performed via an EM algorithm. 
As in the univariate case, we assume that $M_j(\, \cdot \, ; \bfalp_j)$ is a parametric function depending on some vector $\bfalp_j$, $j = 1, \dots, d$, and let $\bfalp = (\bfalp_j)_{j \geq 1}$.
Consider a sample with $N$ clusters of size $d_n$, $n = 1, \dots, N$. Then, the complete log-likelihood is given by  
\begin{align*}
	l _c(\bfalp, \bfpi, \bfT ; (\tilde{\bfy}, \tilde{\bfdelta}) )  & = \sum_{n=1}^{N}  \sum_{j=1}^{d_n}  \Big[\delta_{j,n}\log(\mu_j(y_{j,n}; \bfalp_j)) - z_n M_j(y_{j,n}; \bfalp_j) 
	\\ & \quad \quad \quad \qquad + \log( f_Z(z_n ; \bfpi, \bfT)) \Big]\,,
\end{align*}
where $(\tilde{\bfy}, \tilde{\bfdelta})$ denotes the whole dataset $(y_{j,n}, \delta_{j,n})$, $j = 1, \dots, d_n$, $n = 1, \dots,N $. 

The main difference with respect to the algorithm in Section~\ref{sec:estimation} is in the E-step. We consider one (generic) group of size $d$. First, let us illustrate the calculations by considering only non-censored observations. In such a case, we have that  

\begin{eqnarray*}
	\lefteqn{\E ^{(k+1)}(Z \mid \bfY = \bfy) }\\
	& &= \int_0^\infty z f_{Z|\bfY}(z|\bfy)dz \\
	& &=  \int_0^\infty z \frac{f_{\bfY|Z}(\bfy|z) f_Z(z)}{f_{\bfY}(\bfy)} dz \\
	& &=  \frac{(d+1)! \prod_{j=1}^d \mu_j(y_j;\bfalp^{(k)}_j)}{ f_{\bfY}(\bfy; \bfalp^{(k)}, \bfpi^{(k)}, \bfT^{(k)})  } \bfpi^{(k)} \left(\sum_{j=1}^d M_j(y_j;\bfalp^{(k)}_j) \mat{I} - \bfT^{(k)} \right)^{-2 - d} \bft^{(k)} \\
	& &=  \frac{(d+1) \bfpi^{(k)} \left(\sum_{j=1}^d M_j(y_j; \bfalp^{(k)}_j) \mat{I} - \bfT^{(k)} \right)^{-2 - d} \bft^{(k)} }{ \bfpi^{(k)} \left(\sum_{j=1}^d M_j(y_j; \bfalp^{(k)}_j) \mat{I} - \bfT^{(k)} \right)^{-1 - d} \bft^{(k)}   } \,.
\end{eqnarray*} 
For the general case, where censored observations are present, note that the likelihood of one observation $(\bfy, \bfdelta )$ is given by
\begin{eqnarray*}
\lefteqn{
	 \int_0^{\infty} \prod_{j = 1}^{d} ( z \mu_j(y_j))^{\delta_j} \exp\left( -z M_j(y_j) \right) f_Z(z) dz} \\ 
	& &= \left( \prod_{j = 1}^{d} (  \mu_j(y_j))^{\delta_j} \right) \int_0^{\infty} z^q \exp\left( -z \sum_{j = 1} ^{d} M_j(y_j) \right) f_Z(z) dz \\
	& & = q ! \left( \prod_{j = 1}^{d} (  \mu_j(y_j))^{\delta_j} \right)  \bfpi \left(\sum_{j = 1} ^{d}M_j(y_j) \mat{I} - \bfT \right)^{-1-q} \bft \,,
\end{eqnarray*}
where $q$ is the number of non-censored observations. This yields, 
\begin{align*}
	\E^{(k+1)} (Z \mid (\bfy, \bfdelta))
	& =  \frac{(q+1) \, \bfpi^{(k)}  \left(\sum_{j=1}^d M_j(y_j;\bfalp^{(k)}_j  ) \mat{I} - \bfT^{(k)}  \right)^{-2 - q} \bft^{(k)} }{ \bfpi^{(k)}  \left(\sum_{j = 1}^{d} M_j(y_j; \bfalp^{(k)}_j ) \mat{I} - \bfT^{(k)}  \right)^{-1-q} \bft^{(k)}  } \,.
\end{align*}
This last expression can then be employed in the M-step to find $\bfalp^{(k+1)} $.
Regarding the logarithmic term, the computation of the conditional expectation is similar, and it is easy to see that
\begin{eqnarray*}
	\lefteqn{\E^{(k+1)} (\log (f_Z(Z ; \bfpi, \bfT) \mid (\bfy, \bfdelta)) }\\
	& &= \int \log (f_Z(z ; \bfpi, \bfT))  \frac{ z^q \exp\left( -z \sum_{j =1}^d M_j(y_j; \bfalp^{(k)}_j) \right) \, \bfpi^{(k)}  \exp(\bfT^{(k)} z ) \bft^{(k)} }{ q! \,\bfpi^{(k)}  \left(\sum_{j = 1}^{d} M_j(y_j; \bfalp^{(k)}_j ) \mat{I} - \bfT^{(k)}  \right)^{-1-q} \bft^{(k)}  } dz\,,
\end{eqnarray*}
which can then be maximized in the same way as in the univariate phase-type frailty model to obtain $(\bfpi^{(k+1)} , \bfT^{(k+1)} )$. The detailed routine can be found in Algorithm~\ref{alg:PHshared}. %{\red 
Note that this routine shares the same computational challenges as Algorithm~\ref{alg:PHfrailty}, as described in Remark~\ref{rm:comp}. Additionally, there is a slight increase in computational cost, which grows with the dimension $d$, due to the estimation of additional parameters associated with multiple baseline hazard functions and the computation of matrix powers.
%}

\begin{algorithm}[]
\caption{EM algorithm for the shared phase-type frailty model}\label{alg:PHshared}
\begin{algorithmic}
\normalsize 
\State \textit{\textbf{Input}: Sample $(y_{j,n}, \delta_{j,n})$, $j = 1,\dots, d_n$, $n= 1, \dots, N$.}\\
\begin{enumerate}[label=\arabic*.]
\item[0.] Initialize with some ``arbitrary" $( \bfalp^{(0)}, \bfpi^{(0)} , \bfT^{(0)})$ and let $k = 0$.
	\item \textit{E-step:}  Compute 
\begin{align*}
	\E^{(k+1)} (Z \mid (\bfy_n,\bfdelta_n)) 
	& =  \frac{(q_n + 1) \bfpi^{(k)} \left(\sum_{j = 1}^{d_n}M_j(y_{j,n}; \bfalp^{(k)}_j) \mat{I} - \bfT^{(k)} \right)^{-2 - q_n} \bft^{(k)} }{  \bfpi^{(k)} \left(\sum_{j = 1}^{d_n}M_j(y_{j,n}; \bfalp^{(k)}_j) \mat{I} - \bfT^{(k)} \right)^{-1- q_n} \bft^{(k)} }  \,,
\end{align*}
where $q_n$ is the number of non-censored observations in group $n$, $n = 1,\dots,N$. Let
\begin{align*}
	 g (z )
	&  = \frac{1}{N} \sum _{n=1}^{N}  \frac{z^{q_n} \exp\left(- z \sum_{j = 1}^{d_n}M_j(y_{j,n}; \bfalp^{(k)}_j  )\right)\bfpi^{(k)} \exp(\bfT^{(k)}z)\bft^{(k)}}{ \bfpi^{(k)} \left(\sum_{j = 1}^{d_n}M_j(y_{j,n}; \bfalp^{(k)}_j) \mat{I} - \bfT^{(k)} \right)^{-1- q_n} \bft^{(k)}} \,.
\end{align*}
	
	\item \textit{M-step: Let} 
	\begin{align*}
	\bfalp^{(k +1)} = \argmax_{\bfalp} \sum _{n=1}^{N} \sum _{j=1}^{d_n}  \Bigg[ \delta_{j,n} \log(\mu_j(y_{j,n}; \bfalp_j)) - \E^{(k+1)}  (Z \mid  (\bfy_n, \bfdelta_n)) M_j(y_{j,n}; \bfalp_j) \Bigg] \,,
\end{align*}
assign $\bfeta:= \bfpi^{(k)}$, $\mat{S}:= \bfT^{(k)}$, $\mat{s}:= \bft^{(k)}$, and iterate as necessary the following routine:

\begin{enumerate}
	\item Calculate
\begin{gather*}
	\hat{\eta}_{i} = \int_{0}^{\infty} \dfrac{\eta_{i} \bfe_{i}^{\top} \exp({\bfS y}) \bfs} { \bfeta \exp({\bfS y}) \bfs} g(y)dy \,,\\[4mm]
	\hat{s}_{il} = \dfrac{ \mathlarger{ \int_{0}^{\infty} \dfrac{s_{il}}{ \bfeta \exp({\bfS y}) \bfs } {\int_{0}^{y}  \bfe_{l}^{\top}   \exp({\bfS (y-u)}) \bfs \bfeta \exp({\bfS u}) \bfe_{i} du} \,g(y)dy}}{\mathlarger{\int_{0}^{\infty} \dfrac{1}{ \bfeta \exp({\bfS y}) \bfs}{\int_{0}^{y} \bfe_{i}^{\top} \exp({ \bfS (y-u)}) \bfs \bfeta \exp({\bfS u}) \bfe_{i} du} \, g(y) dy}} \,,\\[4mm]
	\hat{s}_{i} = \dfrac{ \mathlarger{\int_{0}^{\infty} {s_{i}} \dfrac {\bfeta \exp({\bfS y}) \bfe_{i}} {\bfeta \exp({\bfS y}) \bfs } \, g(y) dy}}{\mathlarger{\int_{0}^{\infty} \dfrac{1}{\bfeta \exp({\bfS y}) \bfs }{\int_{0}^{y} \bfe_{i}^{\top} \exp({\bfS (y-u)}) \bfs \bfeta \exp({\bfS u}) \bfe_{i} du} \, g(y) dy}}\,,\\[4mm]
	\hat{s}_{ii} = - \sum_{l\neq i} \hat{s}_{il} - \hat{s}_i \,.
\end{gather*}

Let $\hat{\bfeta} = (\hat{\eta}_{1}, \dots, \hat{\eta}_{p})$, $\hat{\bfS} = \{\hat{s}_{il}\}_{i,l = 1,\dots,p}$, and $\hat{\bfs} = (\hat{s}_1, \dots, \hat{s}_p)^{\top}$.
\item Assign $\bfeta:= \hat{\bfeta}$, $\mat{S}:= \hat{\mat{S}}$, $\mat{s}:= \hat{\mat{s}}$  and GOTO (a).
\end{enumerate}\

	\item Assing  $ \bfpi^{(k + 1)}:=\bfeta$, $\bfT^{(k+1)}:= \mat{S} $, $\bft^{(k+1)}:= \mat{s}$, $k := k+1$, and GOTO 1 until a stopping rule is satisfied. 
\end{enumerate}
\State \textit{\textbf{Output}: Fitted parameters $(\bfalp^{(k)}, \bfpi^{(k)} , \bfT^{(k)})$.}
\end{algorithmic}
\end{algorithm}

\subsection{Correlated phase-type frailty}
Next, we show how correlated frailty can be performed using multivariate phase-type distributions with a focus on the bivariate case for clarity of the presentation. Recall that the correlated frailty model is, on the one hand, an extension of the shared frailty model and, on the other, of the univariate frailty model, where the frailties of individuals in a group are correlated but not necessarily shared. Under this setting, the conditional joint survival function of $\bfY = (Y_1, Y_2)$ given the correlated frailty $\bfZ = (Z_1, Z_2)$ is assumed to be by 
\begin{align*}
	S_{\bfY | \bfZ}(y_1, y_2 \mid Z_1, Z_2 ) & = S_1(y_1 \mid Z_1 ) S_2( y_2 \mid  Z_2 ) \\
	& = \exp(-Z_1 M_1(y_1)) \exp(-Z_2 M_2(y_2)) \,,
\end{align*}
where $M_1$ and $M_2$  denote baseline cumulative hazards. Here, the distribution of the random vector $\bfZ$ determines the association structure of the events in the model. For example, in the particular case $Z_1 = Z_2$, one recovers the shared frailty model. Another notable instance is when $Z_1$ and $Z_2$ are independent, which yields independence of $Y_1$ and $Y_2$.
% independent following univariate frailty specifications. 
Different distributional assumptions for $\bfZ$ have been introduced in the literature. For instance, the correlated Gamma frailty introduced by \cite{yashin1995correlated} possesses closed-form expressions for different functionals. However, it only allows for a positive correlation of the correlated frailties. Another model that allows for a broader range of correlations is the correlated lognormal frailty first introduced by \cite{xue1996bivariate}. Nevertheless, this latter model does not allow for explicit representation of several functionals. The main advantage of the use of multivariate phase-type distributions as frailties is that we obtain closed-form expressions while enabling us to describe any dependence structure, as we discuss next.

 For illustration purposes, we consider the bivariate class of phase-type distributions for which the joint density is given by
\begin{align}\label{eq:bivPH}
	f_{\bfZ}(z_1, z_2) = \bfeta \exp({\bfT_{11} z_1})\bfT_{12}\exp({\bfT_{22} z_2})(-\bfT_{22}) \bfe \,, 
\end{align}
where $\bfeta$ is a $p_1$-dimensional vector of initial probabilities, $\bfT_{11}$ and $\bfT_{22}$ are sub-intensity matrices of dimensions $p_1 \times p_1$ and $p_2 \times p_2$, respectively, and $\bfT_{12} $ is a $p_1 \times p_2$ non-negative matrix satisfying $\bfT_{11} \bfe  = -\bfT_{12} \bfe $. In particular, this explicit form of the joint density yields the following parametrization of the phase-type margins: $Z_1 \sim  \mbox{PH}(\bfeta , \bfT_{11} )$ and $Z_2 \sim  \mbox{PH}(\bfeta(\bfT_{12})^{-1} , \bfT_{22} )$. Moreover, it is essential to note that the case $Z_1 \sim  \mbox{PH}(\bfpi , \bfT )$ and $Z_2 \sim  \mbox{PH}(\vect{\nu}, \bfS )$ independent is retrieved with $\bfeta  = (\bfpi, \0)$, $\bfT_{11} = \bfT$, $\bfT_{22} = \bfS$, and $\bfT_{12} = \bft \vect{\nu}$. For our purposes, the following expression derived straightforwardly from \eqref{eq:bivPH} is integral
\begin{align*}
	 \E( \exp({-u_1 Z_1 - u_2 Z_2})) =\bfeta ( u_1 \mat{I} - \bfT_{11})^{-1}\bfT_{12}( u_2 \mat{I} - \bfT_{22})^{-1}(-\bfT_{22}) \bfe  \,.
\end{align*}

In addition to closed-form expressions of different functionals, this class of bivariate distributions is dense on the set of distributions on $\R_+^2$ \cite[see][]{albrecher2020fitting}. Thus, we can approximate any distribution with support on  $\R_+^2$ with members of this family. Moreover, these distributions allow for estimation via an explicit EM algorithm introduced by \cite{ahlstrom1999parametric}, which facilitates their use in practice. It is worth mentioning that this class can be extended to larger dimensions with a similar construction principle, and we refer to, e.g., \cite{albrecher2020mulmml} for the interested reader. However, to focus on the main ideas of the model, we consider only the bivariate case. 

We now take $ \bfZ$ with joint density given by \eqref{eq:bivPH}. In this way, $\bfY$ has joint survival function given by
\begin{align*}
	S_{\bfY}(y_1, y_2 ) &= \int_0^\infty \int_0^\infty  \exp(-z_1 M_1(y_1)) \exp(-z_2 M_2(y_2))  f_{\bfZ}(z_1, z_2) dz_1 dz_2 \\
	&= \bfeta (  M_1(y_1) \mat{I} - \bfT_{11})^{-1}\bfT_{12}(  M_2(y_2) \mat{I} - \bfT_{22})^{-1}(-\bfT_{22}) \bfe  \,.
\end{align*}

Moreover, from this last expression, it follows that 
\begin{align*}
	f_{\bfY}(y_1, y_2 ) = \mu_1(y_1) \mu_2(y_2) \bfeta (  M_1(y_1) \mat{I} - \bfT_{11})^{-2}\bfT_{12}(  M_2(y_2) \mat{I} - \bfT_{22})^{-2}(-\bfT_{22}) \bfe  \,.
\end{align*}

\subsubsection{Estimation} Estimation of this model can also be carried out via an EM algorithm. 
To keep the presentation concise, we present the calculations only for the case of non-censored observations. Nevertheless,  the case of right-censoring follows similarly. Again, we assume that $M_j(\cdot ; \bfalp_j)$ is a parametric function depending on some vector $\bfalp_j$, $j=1,2$, and let $\bfalp = (\bfalp_1, \bfalp_2)$. Then, for a sample of size $N$ with observations $(y_{1,1},y_{2,1}), \dots, (y_{1,N},y_{2,N})$, the complete log-likelihood is given by 
\begin{eqnarray*}
	 \lefteqn{l_c(\bfalp, \bfeta, \bfT_{11} , \bfT_{12},  \bfT_{22} ; \tilde{\bfy} )}\\
	  	& & = \sum_{n= 1}^{N} \Big[  \log( \mu_1(y_{1,n}; \bfalp_1)) + \log( \mu_2(y_{2,n}; \bfalp_2)) -z_{1,n} M_1(y_{1,n}; \bfalp_1)  \\
	& & \quad \quad \qquad -z_{2,n} M_2(y_{2,n}; \bfalp_2) + \log( f_{\bfZ}(z_{1,n}, z_{2,n} ; \bfeta, \bfT_{11} , \bfT_{12},  \bfT_{22} ) ) \Big]\,.
\end{eqnarray*}
First, we compute the conditional expectations given the observed data by considering one (generic) data point $\bfy =(y_1, y_2)$. Note that
\begin{eqnarray*}
	\lefteqn{\E ^{(k+1)} (Z_1 \mid \bfY = \bfy) } \\
	& &= \int_0^\infty \int_0^\infty  z_1 \frac{f_{\bfY | \bfZ}(y_1, y_2 | z_1, z_2) }{f_{\bfY}(y_1, y_2 )} f_{\bfZ}(z_1, z_2) dz_1 dz_2 \\
	& &= \frac{\mu_1(y_{1}) \mu_2(y_{2})}{f_{\bfY}(y_1, y_2 )}\int_0^\infty \int_0^\infty  z_1^2 z_2 \exp(-z_1 M_1(y_{1})) \exp(-z_2 M_2(y_{2})) f_{\bfZ}(z_1, z_2) dz_1 dz_2 \\
	& &= \frac{\mu_1(y_{1} ) \mu_2(y_{2})}{f_{\bfY}(y_1, y_2 )} 2\bfeta^{(k)}  (  M_1(y_1) \mat{I} - \bfT_{11}^{(k)} )^{-3}\bfT_{12}^{(k)} (  M_2(y_2) \mat{I} - \bfT_{22}^{(k)} )^{-2}(-\bfT_{22}^{(k)} ) \bfe \\
	& &= \frac{2 \bfeta^{(k)}  (  M_1(y_1; \bfalp^{(k)}_1) \mat{I} - \bfT_{11}^{(k)} )^{-3}\bfT_{12}^{(k)} (  M_2(y_2; \bfalp^{(k)}_2) \mat{I} - \bfT_{22}^{(k)} )^{-2}(-\bfT_{22}^{(k)} ) \bfe}{\bfeta^{(k)}  (  M_1(y_1; \bfalp^{(k)}_1) \mat{I} - \bfT_{11}^{(k)} )^{-2}\bfT_{12}^{(k)} (  M_2(y_2; \bfalp^{(k)}_2) \mat{I} - \bfT_{22}^{(k)} )^{-2}(-\bfT_{22}^{(k)} ) \bfe} \,.
\end{eqnarray*}
Similarly,
\begin{eqnarray*}
	\lefteqn{\E^{(k+1)}(Z_2 \mid \bfY = \bfy) }\\
	& & = \frac{2 \bfeta^{(k)} (  M_1(y_1; \bfalp^{(k)}_1) \mat{I} - \bfT_{11}^{(k)})^{-2}\bfT_{12}^{(k)}(  M_2(y_2; \bfalp^{(k)}_2) \mat{I} - \bfT_{22}^{(k)})^{-3}(-\bfT_{22}^{(k)}) \bfe}{\bfeta^{(k)} (  M_1(y_1; \bfalp^{(k)}_1) \mat{I} - \bfT_{11}^{(k)})^{-2}\bfT_{12}^{(k)}(  M_2(y_2; \bfalp^{(k)}_2) \mat{I} - \bfT_{22}^{(k)})^{-2}(-\bfT_{22}^{(k)}) \bfe} \,.
\end{eqnarray*}
For the logarithmic term, we have that 
\begin{eqnarray*}
	\lefteqn{\E (\log( f_{\bfZ}(Z_1, Z_2 ; \bfeta, \bfT_{11} , \bfT_{12},  \bfT_{22})) \mid \bfY = \bfy) }\\
	   & & = \int_0^\infty \int_0^\infty  \log( f_{\bfZ}(z_1, z_2 ; \bfeta, \bfT_{11} , \bfT_{12},  \bfT_{22})) f_{\bfZ | \bfY}(z_1, z_2 | y_1, y_2)  dz_1 dz_2 \\
	    & & = \int_0^\infty \int_0^\infty  \log( f_{\bfZ}(z_1, z_2 ; \bfeta, \bfT_{11} , \bfT_{12},  \bfT_{22}))  z_1 z_2 \exp\left(-\sum_{j =1}^2z_j M_j(y_j;\bfalp^{(k)}_j ) \right)  \times\\
	    & & \quad \frac{ \bfeta^{(k)} \exp({\bfT_{11}^{(k)} z_1})\bfT_{12}^{(k)}\exp({\bfT_{22}^{(k)} z_2})(-\bfT_{22}^{(k)}) \bfe}{ \bfeta^{(k)} (  M_1(y_1;\bfalp^{(k)}_1) \mat{I} - \bfT_{11}^{(k)})^{-2}\bfT_{12}^{(k)}(  M_2(y_2;\bfalp^{(k)}_2) \mat{I} - \bfT_{22}^{(k)})^{-2}(-\bfT_{22}^{(k)}) \bfe}  dz_1 dz_2 \,.
\end{eqnarray*} 
Once the conditional expectations are computed, $\bfalp^{(k+1)}$ can be found numerically as in the univariate phase-type frailty model. To obtain the updated parameters $\bfeta^{(k+1)}$, $\bfT_{11}^{(k+1)}$, $\bfT_{12}^{(k+1)}$,  and $\bfT_{22}^{(k+1)}$ of the bivariate phase-type component, the maximization needed can be performed by fitting a bivariate phase-type distribution to an appropriate joint density. In \cite{albrecher2020fitting}, it is shown how the EM algorithm in \cite{ahlstrom1999parametric} can be modified for such a purpose.  We have delineated this procedure in Algorithm~\ref{alg:bivPH} for the sake of completeness. 
Finally, we detail the steps of the EM algorithm for the correlated phase-type frailty model in Algorithm~\ref{alg:PHcor}, which includes the case of right-censored observations.
%{\red
For a detailed application of this model in the context of mixed Poisson regression, we refer to \cite{furrer2024bivariate}.
%, where the authors demonstrate its superior modeling performance over other traditional specifications.
%}

\begin{algorithm}[]
\caption{EM algorithm for the correlated phase-type frailty model}\label{alg:PHcor}
\begin{algorithmic}
\normalsize 
\State \textit{\textbf{Input}: Sample $(y_{j,n}, \delta_{j,n})$, $j = 1, 2$, $n= 1, \dots, N$.}\\
\begin{enumerate}[label=\arabic*.]
\item[0.] Initialize with some ``arbitrary" $( \bfalp^{(0)}, \bfeta^{(0)} , \bfT^{(0)}_{11}, \bfT^{(0)}_{12}, \bfT^{(0)}_{22})$ and let $k = 0$.
	\item \textit{E-step:}  Compute 
\begin{eqnarray*}
	\lefteqn{\E^{(k+1)}(Z_1 \mid (\bfy_n, \bfdelta_n)) }\\
	& & = \frac{(1 + \delta_{1,n}) \bfeta^{(k)} (  M_1(y_{1,n}; \bfalp^{(k)}_1) \mat{I} - \bfT_{11}^{(k)})^{-2-\delta_{1,n}}\bfT_{12}^{(k)}(  M_2(y_{2,n}; \bfalp^{(k)}_2) \mat{I} - \bfT_{22}^{(k)})^{-1-\delta_{2,n}}(-\bfT_{22}^{(k)}) \bfe}{\bfeta^{(k)} (  M_1(y_{1,n}; \bfalp^{(k)}_1) \mat{I} - \bfT_{11}^{(k)})^{-1-\delta_{1,n}}\bfT_{12}^{(k)}(  M_2(y_{2,n}; \bfalp^{(k)}_2) \mat{I} - \bfT_{22}^{(k)})^{-1-\delta_{2,n}}(-\bfT_{22}^{(k)}) \bfe} \,,
\end{eqnarray*}

\begin{eqnarray*}
	\lefteqn{\E^{(k+1)}(Z_2 \mid (\bfy_n, \bfdelta_n)) }\\
	& & = \frac{(1 + \delta_{2,n}) \bfeta^{(k)} (  M_1(y_{1,n}; \bfalp^{(k)}_1) \mat{I} - \bfT_{11}^{(k)})^{-1 - \delta_{1,n}}\bfT_{12}^{(k)}(  M_2(y_{2,n}; \bfalp^{(k)}_2) \mat{I} - \bfT_{22}^{(k)})^{-2-\delta_{2,n}}(-\bfT_{22}^{(k)}) \bfe}{\bfeta^{(k)} (  M_1(y_{1,n}; \bfalp^{(k)}_1) \mat{I} - \bfT_{11}^{(k)})^{-1-\delta_{1,n}}\bfT_{12}^{(k)}(  M_2(y_{2,n}; \bfalp^{(k)}_2) \mat{I} - \bfT_{22}^{(k)})^{-1 - \delta_{2,n}}(-\bfT_{22}^{(k)}) \bfe} \,,
\end{eqnarray*}
and let
\begin{eqnarray*}
	\lefteqn{ g (\bfz )}\\
	& &  = \frac{1}{N} \sum _{n=1}^{N} \frac{ z_1^{\delta_{1,n}} z_2^{\delta_{2,n}} \exp\left(-\sum_{j =1}^2z_j M_j(y_{j,n};\bfalp^{(k)}_j ) \right)  \bfeta^{(k)} \exp({\bfT_{11}^{(k)} z_1})\bfT_{12}^{(k)}\exp({\bfT_{22}^{(k)} z_2})(-\bfT_{22}^{(k)}) \bfe}{ \bfeta^{(k)} (  M_1(y_{1,n};\bfalp^{(k)}_1) \mat{I} - \bfT_{11}^{(k)})^{-1 - \delta_{1,n}}\bfT_{12}^{(k)}(  M_2(y_{2,n};\bfalp^{(k)}_2) \mat{I} - \bfT_{22}^{(k)})^{-1-\delta_{2,n}}(-\bfT_{22}^{(k)}) \bfe}  \,.
\end{eqnarray*}
	
	\item \textit{M-step: Let} 
	\begin{align*}
	\bfalp^{(k +1)} = \argmax_{\bfalp} \sum _{n=1}^{N}  \sum _{j=1}^{2} \Big[ \delta_{j,n}\log(\mu_j(y_{j,n}; \bfalp_j)) - \E^{(k+1)}  (Z_j \mid (\bfy_n, \bfdelta_n))  M_j(y_{j,n}; \bfalp_j) \Big] \,,
\end{align*}
assign $\bfeta:= \bfeta^{(k)}$, $\mat{T}_{11}:= \bfT^{(k)}_{11}$, $\mat{T}_{12}:= \bfT^{(k)}_{12}$, $\mat{T}_{22}:= \bfT^{(k)}_{22}$, and iterate as necessary the  routine detailed in Algorithm~\ref{alg:bivPH}.

	\item Assing  $ \bfeta^{(k + 1)}:=\bfeta$, $\bfT^{(k+1)}_{11}:= \mat{T}_{11} $, $\bfT^{(k+1)}_{12}:= \mat{T}_{12} $, $\bfT^{(k+1)}_{22}:= \mat{T}_{22} $, $k := k+1$, and GOTO 1 until a stopping rule is satisfied. 
\end{enumerate}
\State \textit{\textbf{Output}: Fitted parameters $( \bfalp^{(k)}, \bfeta^{(k)} , \bfT^{(k)}_{11}, \bfT^{(k)}_{12}, \bfT^{(k)}_{22})$.}
\end{algorithmic}
\end{algorithm}

%{\red 

\begin{remark}[On higher-dimensional extensions and computational challenges] \rm 
		We have presented the theoretical foundations and derived a fitting algorithm for the correlated phase-type frailty model in the bivariate case when considering an underlying multivariate phase-type model with density \eqref{eq:bivPH}. However, a similar construction can be carried out using other families of multivariate phase-type distribution. In particular, one can employ the class of mPH distributions introduced in \cite{bladt2022mph}, which is, in fact, a subfamily of the distributions described by \eqref{eq:bivPH}. An advantage of considering the mPH specification is that the theoretical developments of the resulting correlated frailty model can be carried out straightforwardly to higher dimensions. Nevertheless, estimating these high-dimensional models using algorithms analogous to  Algorithm~\ref{alg:PHcor} would be highly computationally demanding. In fact, Algorithm~\ref{alg:PHcor} is the most computationally intensive among the routines derived here. This complexity arises from two main sources. First, the model involves more parameters than the other specifications, and it requires the computation of matrix exponentials for two different sub-intensity matrices. Second, applying Algorithm~\ref{alg:bivPH} to fit the joint density $g$, as described in Algorithm~\ref{alg:PHcor}, introduces an additional computational burden due to the presence of multiple integrals, a type of computation that suffers from the curse of dimensionality.
\end{remark}

\section{Examples}\label{sec:examples}
This section presents several detailed numerical illustrations of the estimation of the phase-type frailty models presented in previous sections. In all the examples, the number of iterations in the algorithms was chosen so that changes in successive log-likelihoods became negligible, which is a typical stopping rule when working with the EM algorithm. 
It is essential to emphasize that the central objective of this section is not to undertake an exhaustive comparison with, or establish superiority over, other models. Instead, our focus remains on highlighting the computational practicability and flexibility of the proposed models. %{\red 
To support the practical application of these models, a preliminary version of an R package implementing the proposed algorithms, along with the code used for the following illustrations, is available in the Zenodo repository \url{https://zenodo.org/records/15045852}.
\subsection{Phase-type-Gompertz frailty model}
\begin{figure}[h]
\centering
\includegraphics[width=0.49\textwidth]{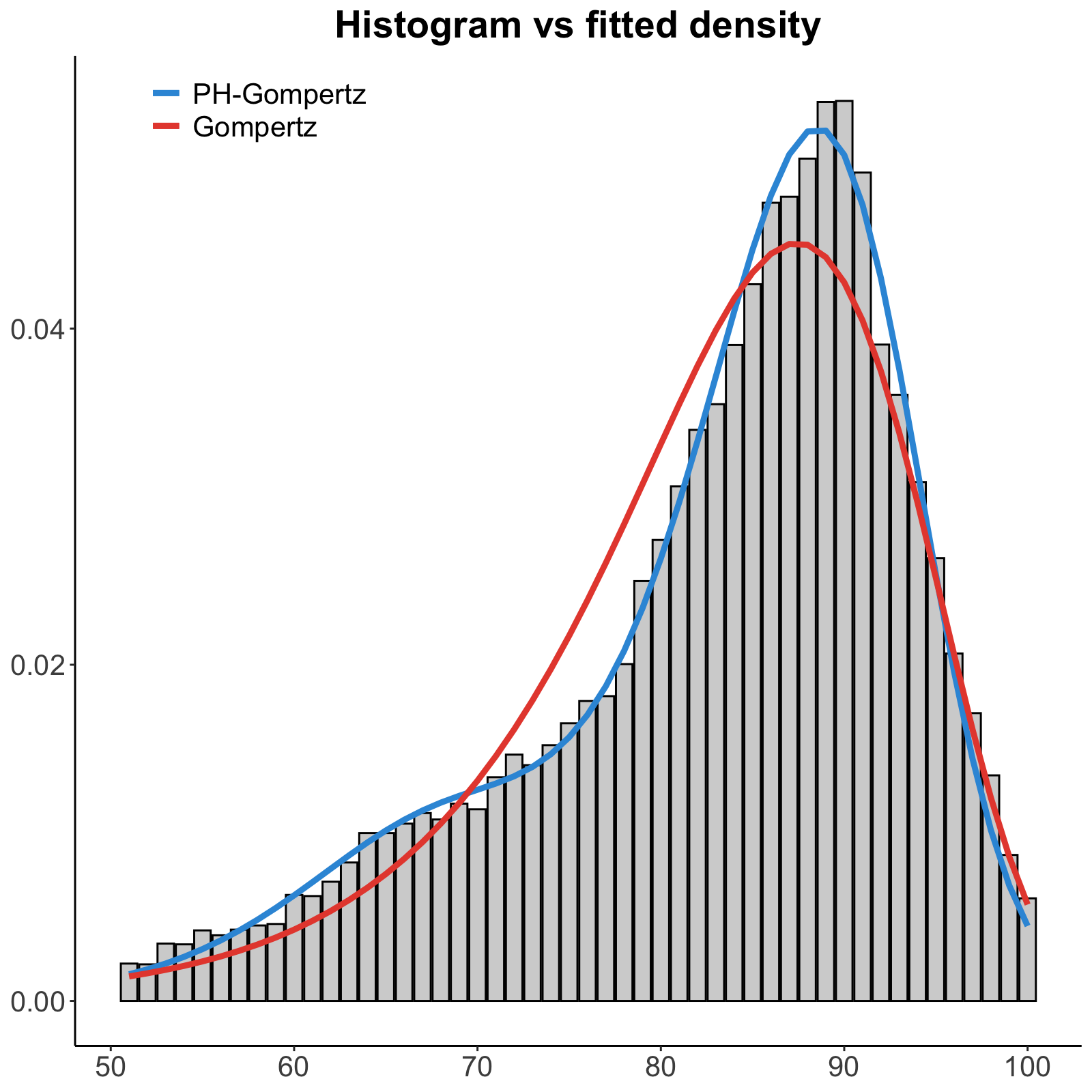}
\includegraphics[width=0.49\textwidth]{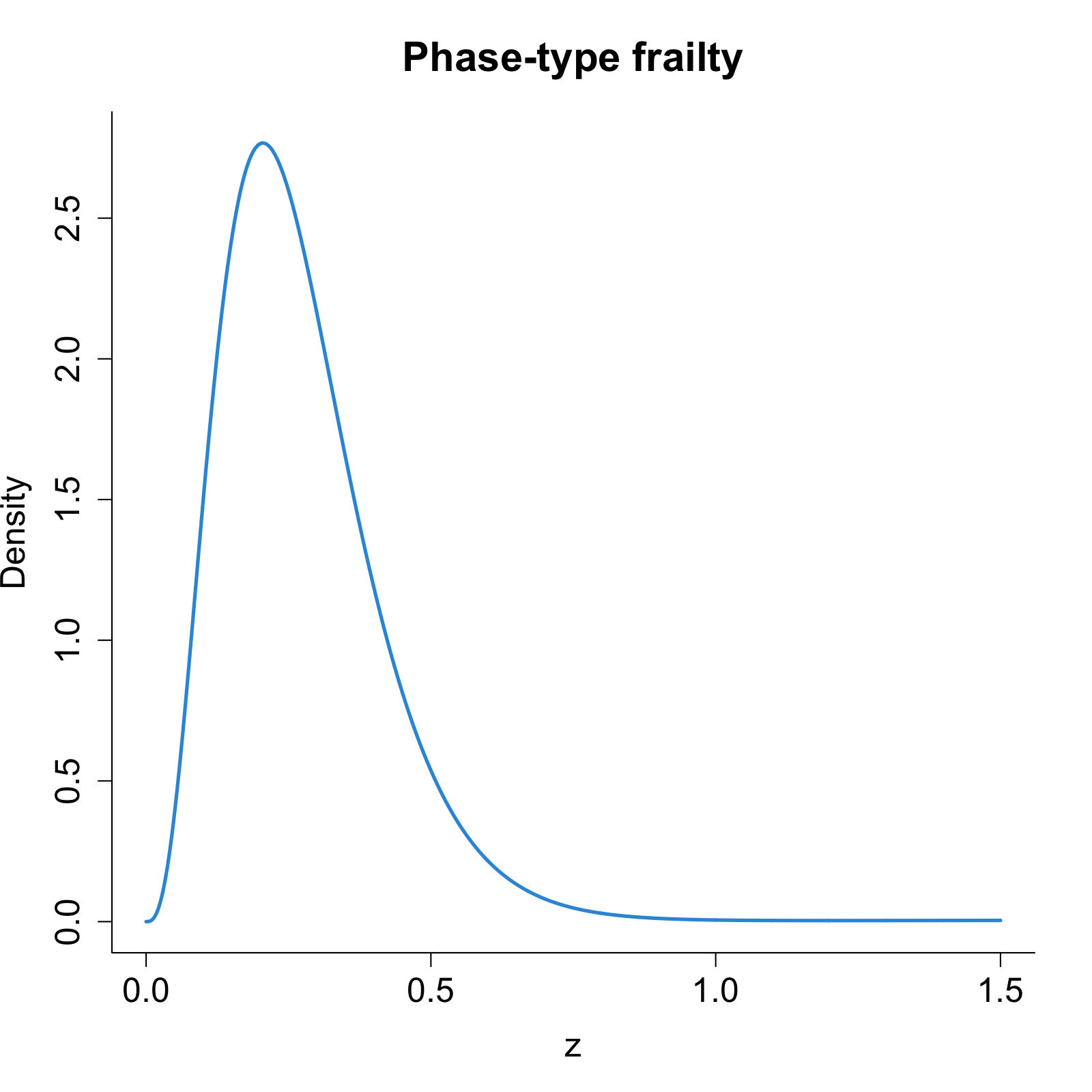}
\caption{Histogram of lifetimes of the Swedish female population that died in the year 2011 at ages 50 to 100 versus density of the fitted phase-type-Gompertz frailty model and density of fitted Gompertz distribution (left) and density function of the underlying phase-type frailty (right).}
\label{fig:phgompertz}
\end{figure}

%%{\red Expand a little bit about the importance of frailty models in mortality modeling}
It is well-known that frailty models with Gompertz baseline hazards, that is $\mu(y) = b \exp({cy})$ with $b,c>0$, are useful tools in describing human mortality at adult and old ages, being the Gamma distribution the most common frailty employed in practice \cite[see, e.g.,][]{butt2004application,missov2013gamma,vaupel1979impact}. 
%%It is well-known that the Gompertz-related parametric models are useful in describing the pattern of adult human deaths. 
%%The gamma-Gompertz multiplicative frailty model is the most common parametric model applied to human mortality data at adult and old ages.
%%Frailty models with baseline Gompertz hazards, $\mu(y) = b \exp({cy})$, have been employed for mortality modeling in the literature, being the most popular choice for the frailty the Gamma distribution. The objective is to account for unobserved heterogeneity in the population.... For instances...
%
In this example, we propose the use of phase-type frailty models with baseline Gompertz hazard to model this type of data. As a specific case of study, we consider the lifetimes of the Swedish female population that died in the year $2011$ at ages $50$-$100$. This data was obtained from the Human Mortality Database (HMD) and available in the R-package \texttt{MortalitySmooth} \citep{camarda2012mortalitysmooth}. Then, we fit a  univariate phase-type frailty model with a Coxian structure of 6 phases for the phase-type component, obtaining the following estimated parameters
%%\begin{gather*} 
%%	\hat{\bfpi}=\left(
%%	1, \,0 , \,0 , \,0 , \,0 , \,0 \right)\,, \\ 
%%	\hat{\bfT}=\left( \begin{array}{cccccc}
%%	-14.6496 &  14.6496 & 0 & 0 & 0 & 0 \\
%%	0 &  -14.6496 & 14.6496 & 0 & 0 & 0 \\
%%	0 &  0 & -14.6496 & 14.6496 & 0 & 0 \\
%%	0 &  0 & 0 & -14.6496 & 2.2944 & 0 \\
%%	0 &  0 & 0 & 0 & -0.4105 & 0.4105 \\
%%	0 &  0 & 0 & 0 & 0 & -0.0774 
%%	\end{array} \right) \,, \\ 
%%	\hat{b}= 0.000497 \quad \hat{c}= 0.187532 \,.
%%\end{gather*}
 \begin{gather*} 
	\hat{\bfpi}=\left(
	1, \,0 , \,0 , \,0 , \,0 , \,0 \right)\,, \\ 
	\hat{\bfT}=\left( \begin{array}{cccccc}
	-14.6288 &  14.6288 & 0 & 0 & 0 & 0 \\
	0 &  -14.6288 & 14.6288 & 0 & 0 & 0 \\
	0 &  0 & -14.6288 & 14.6288 & 0 & 0 \\
	0 &  0 & 0 & -14.6288 & 2.2798 & 0 \\
	0 &  0 & 0 & 0 & -0.0786 & 0.0786 \\
	0 &  0 & 0 & 0 & 0 & -0.3820 
	\end{array} \right) \,, \\ 
	\hat{b}= 0.000499\,, \quad \hat{c}= 0.187364 \,.
\end{gather*}
Figure~\ref{fig:phgompertz} shows that the fitted distribution provides an adequate model for the sample with better performance than a conventional Gompertz distribution. This is also supported by a comparison of the log-likelihoods of both models, $-162,403.6$ for the Gompertz distribution and $-161,769.9$ for the phase-type-Gompertz model. %{\red 
Additionally, we fitted frailty models with Gamma and inverse Gaussian frailties. However, they show only marginal improvements in fit compared to the standalone Gompertz model, with corresponding log-likelihood values of $-16,2392.2$ and $-162,392.3$, respectively. Moreover, the fitted densities are nearly indistinguishable from that of the Gompertz model and hence omitted in Figure~\ref{fig:phgompertz}.

%}

\subsection{Loss insurance data} 
This example shows that the phase-type frailty model leads to parsimonious distributions that can be employed in contexts beyond survival analysis. Specifically, consider $\mu(y) = \theta y^{\theta-1}$, $\theta>0$, which yields 
\begin{align*}
	S_Y(y) &=  \bfpi ( y^{\theta} \mat{I} - \bfT )^{-1} \bft  = \bfpi ( y^{\theta}  (-\bfT)^{-1}  + \mat{I})^{-1} \bfe \,.
\end{align*}
We call this the {\em matrix-Pareto type III} distribution due to its similarity with the conventional Pareto type III distribution and also to distinguish it from the matrix-Pareto models introduced in \cite{albrecher2019inhomogeneous} and \cite{albrecher2021cph}. 
 For this case, \eqref{eq:asymlaplace} implies that
\begin{align*}
	S_Y(y) \sim D y^{-m \theta } \,, \quad y \to \infty \,,
\end{align*}
meaning that the distribution is regularly varying with index $m\theta$. 
As noted earlier in \cite{albrecher2021cph} for the matrix-Pareto type I and II models,  a primary difference among the matrix-Pareto models is how the tail behavior is specified. 
For the matrix-Pareto type I in \cite{albrecher2019inhomogeneous}, the tail behavior is determined by the largest real eigenvalue of $\bfT$. In contrast, for the matrix-Pareto type II in \cite{albrecher2021cph}, the tail behavior is specified by a scalar (shape) parameter. Finally, for the matrix-Pareto type III introduced here, the tail behavior depends on a scalar parameter and the  sizes of the Jordan blocks of $\bfT$. 
Another critical difference arises when working with covariate information. In the matrix-Pareto type III, the covariates can be incorporated to act multiplicatively in scale. The same applies to the matrix-Pareto type II through the regression setting in 
\cite{bladt2022heavy}. In contrast, the matrix-Pareto type I distribution can include covariates affecting multiplicatively in shape (tail) via the survival regression model in \cite{bladt2020survival}.

\begin{figure}[h]
\centering
\includegraphics[width=0.49\textwidth]{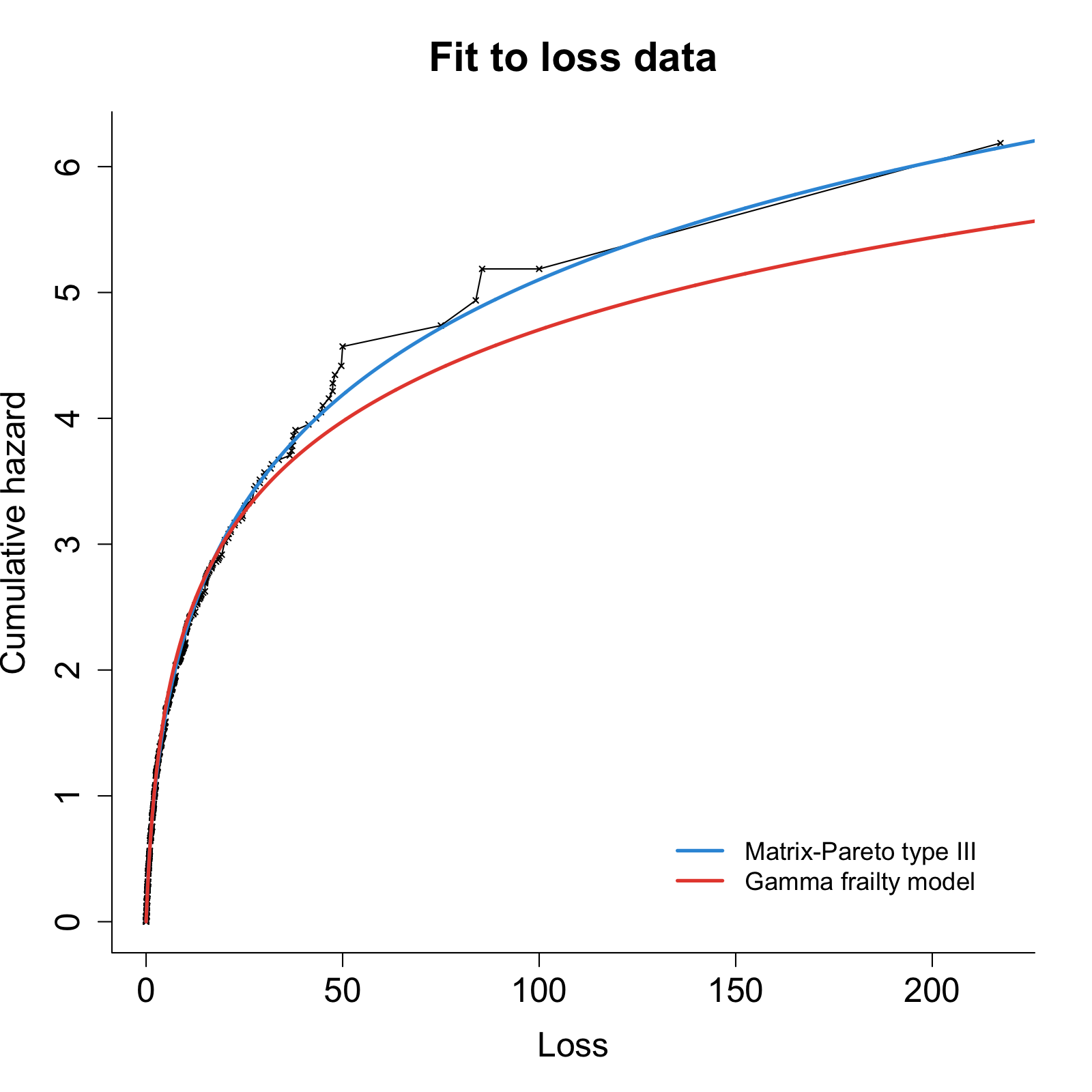}
\includegraphics[width=0.49\textwidth]{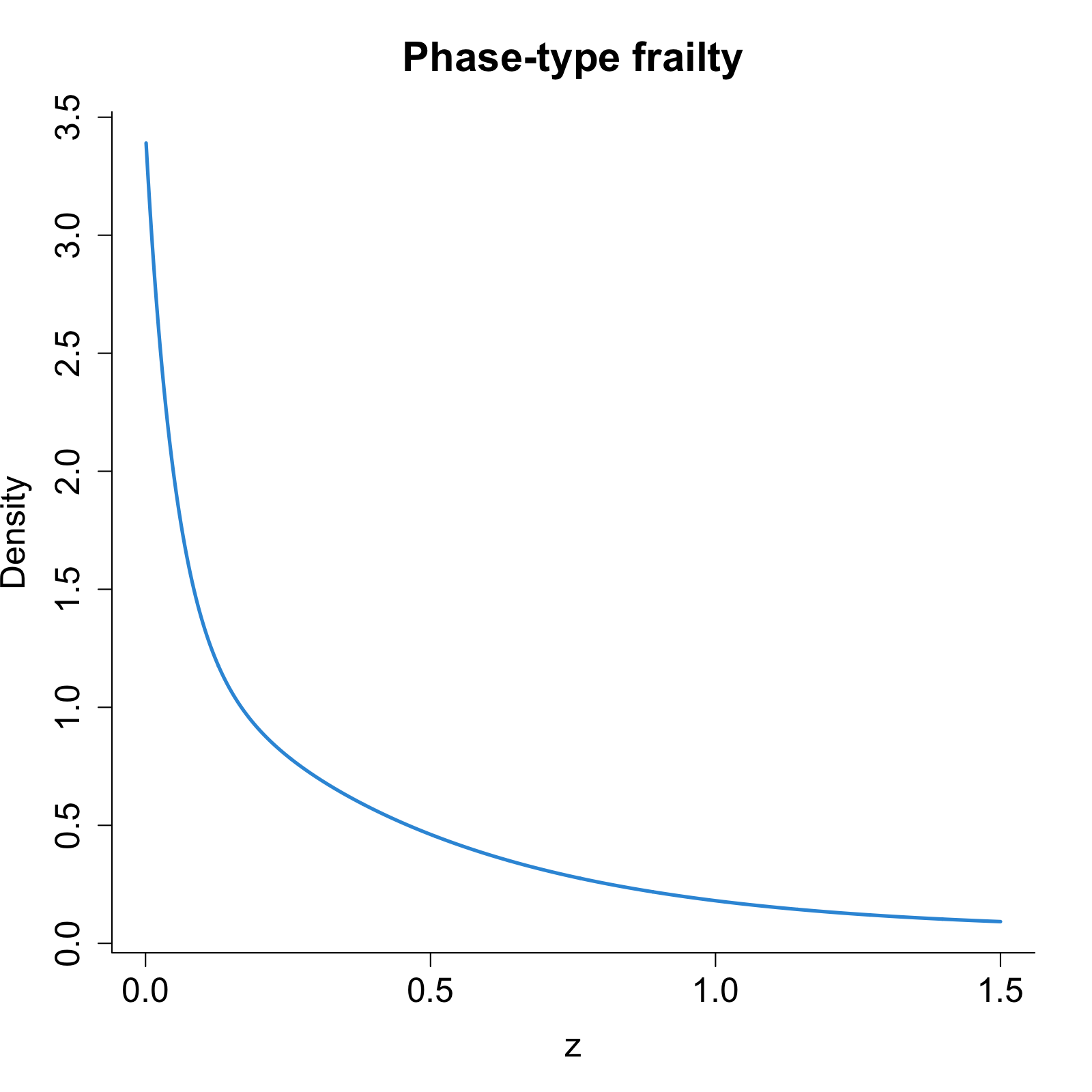}
\caption{Cumulative hazard functions of the fitted matrix-Pareto type III distribution and fitted Gamma frailty model versus the non-parametric Nelson-Aalen estimator of the sample (left) and density function of the underlying phase-type frailty (right) }
\label{fig:loss}
\end{figure}

We now illustrate this model's use in a real-life dataset, namely the loss insurance claim data available in the R-package \texttt{copula} \citep{copularpackage}. The data consists of $1500$ insurance claims from a real-life insurance company, where each data point is conformed of an indemnity payment (\textit{loss}) and an allocated loss adjustment expense (\textit{alae}). For this analysis, we consider only the loss component (scaled by a factor of $10^{-4}$), for which $34$ observations are right-censored. Then, we fit a matrix-Pareto type III distribution to the resulting sample. To reduce the number of parameters, we consider a Coxian structure of dimension $4$ in the phase-type frailty, obtaining in this way the following  estimated parameters
\begin{gather*} 
	\hat{\bfpi}=\left(
	1, \,0,\, 0 \,, 0\right)\,, \\ 
	\hat{\bfT}=\left( \begin{array}{cccc}
	-19.7212 &  16.2868 & 0 & 0 \\
	0 & -2.1507 & 0.7222 & 0  \\
	0 & 0 & -0.5009 & 0.5008 \\
	0 & 0 & 0& -0.5009
	\end{array} \right) \,, \\ 
	\hat{\theta}= 1.3709 \,,
\end{gather*}
%{\red 
For comparison purposes, we also consider a Gamma frailty model with Weibull baseline hazard.
Figure~\ref{fig:loss}  shows the cumulative hazard functions of the two fitted models, where we observe that the matrix-Pareto type III model is closer to the non-parametric Nelson-Aalen estimator, indicating a better fit to the data. The values of the log-likelihoods  also support this conclusion with -3,027.2 for the matrix-Pareto type III model and -3,034.3 for the Gamma frailty model.  %Figure~\ref{fig:loss}  shows that the cumulative hazard of the matrix-Pareto type III model is close to the non-parametric Nelson-Aalen estimator, indicating that the fitted distribution provides an adequate model for the sample. 
Finally, %}
we would also like to mention that an analysis of the same data set employing the matrix-Pareto type I distribution can be found in \cite{bladt2021matrixdist}, and an analysis with the matrix-Pareto type II was done in \cite{albrecher2021cph}. 

\subsection{Lognormal shared frailty}
 The lognormal distribution is among the most important and often-used frailty distributions in practice. In particular, the correlated lognormal frailty model is a relevant tool for modeling dependence in multivariate settings due to its relationship with the multivariate normal distribution. Nevertheless, the univariate lognormal distribution has also been employed in the shared and univariate settings \citep[cf., e.g.,][]{flinn1982new,mcgilchrist1991regression}. However, one downside of using the lognormal distribution as frailty is that there is no closed-form solution for its Laplace transform, and consequently, neither for the survival and density functions of  the frailty model. Hence, numerical approximations of this Laplace transform must be employed (see, e.g., \cite{asmussen2016laplace} for an approach in terms of the Lambert W function). 
 In the following, we illustrate the practical feasibility of the shared phase-type frailty model by fitting it to a simulated bivariate dataset from the shared lognormal frailty model. 
 %In the following, we illustrate the practical feasibility of the shared phase-type frailty model by fitting it to a simulated bivariate dataset from a shared lognormal frailty model. 

%We start with the shared frailty specification. 
 Assume a simple two-group setup of pairs of individuals whose lifetimes follow the following shared lognormal frailty model: the baseline hazards of the individuals are 
Gompertz hazards $\mu_j(y) = b_j \exp({c_jy})$, $j = 1,2$, with $b_1 = 0.01$ and $c_1 = 1$ for the first individual (Individual 1), and  $b_2 = 0.1$ and $c_2 = 2$ for the second individual (Individual 2). The dependence of their lifetimes is modeled with a shared lognormal frailty  with density 
\begin{align*} 
	f_Z(z) = \frac{1}{z \sigma \sqrt{2 \pi}} \exp \left( - \frac{(\log(z) - \nu )^2}{2 \sigma^2} \right) \,, \quad z>0 \,,
\end{align*}
and parameters $\nu = -0.35$ and $\sigma = 0.8$.  
 The first group (Group 1) follows exactly the specification above, while the second group (Group 2) is assumed to have larger proportional hazards by a factor of $ \exp(0.5)$. 
If we consider the categorical covariate $X \in \{ 0,1\}$ and take $\beta = 0.5$, we can summarize the model as follows:
 
Group 1 $(X=0)$
 \begin{align*}
	\text{Individual 1}: \mu_1(t; Z, X) = Z \mu_1(t) \,,\\ 
	\text{Individual 2}: \mu_2(t; Z, X) = Z \mu_2(t) \,.
\end{align*}

Group 2 $(X=1)$
 \begin{align*}
	\text{Individual 1}: \mu_1(t; Z, X) = Z \mu_1(t) \exp(0.5) \,,\\ 
	\text{Individual 2}: \mu_2(t; Z, X) = Z \mu_2(t) \exp(0.5) \,. 
\end{align*}
Next, we simulate $1,000$  bivariate observations for each group. 

\begin{figure}[h]
\centering
\includegraphics[width=0.49\textwidth]{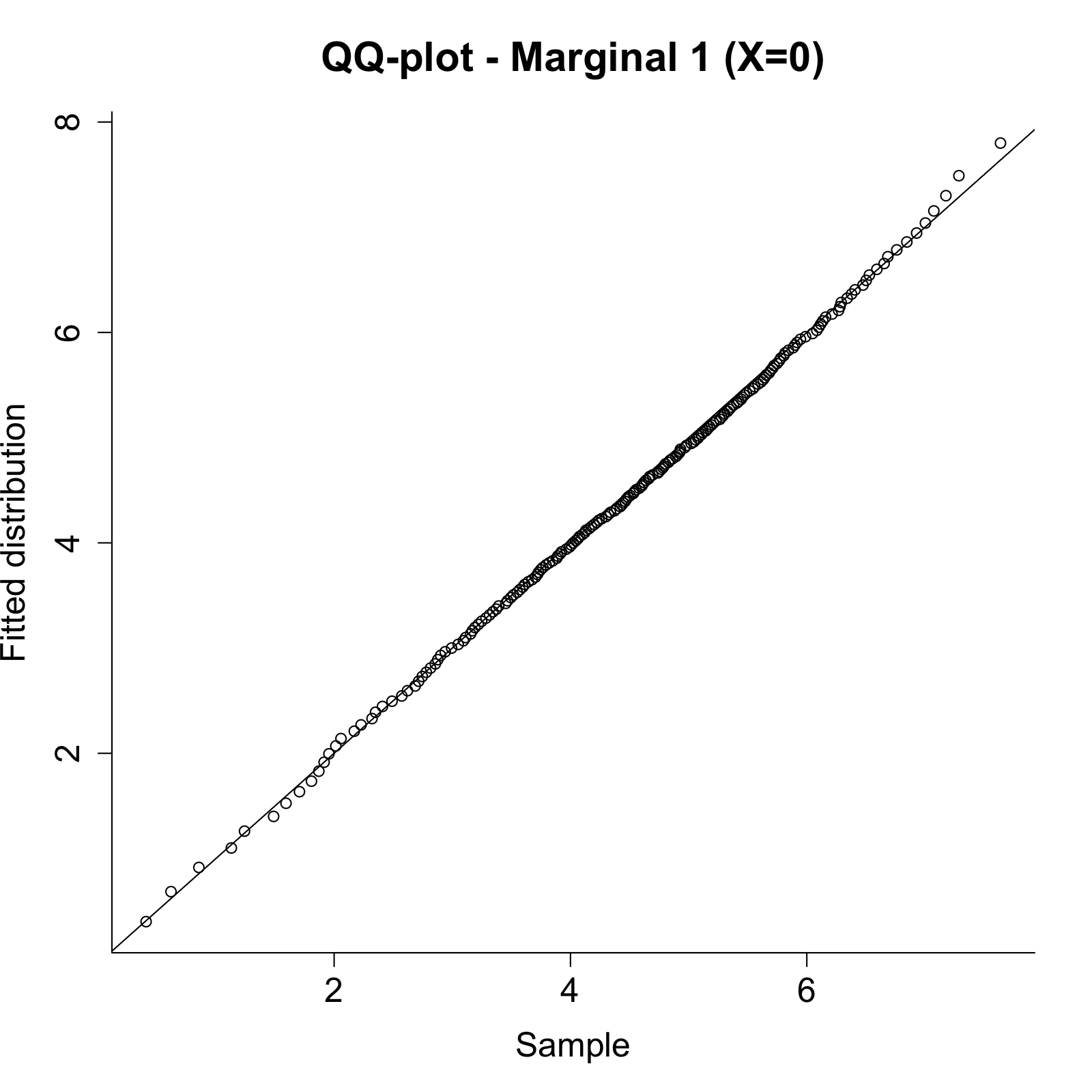}
\includegraphics[width=0.49\textwidth]{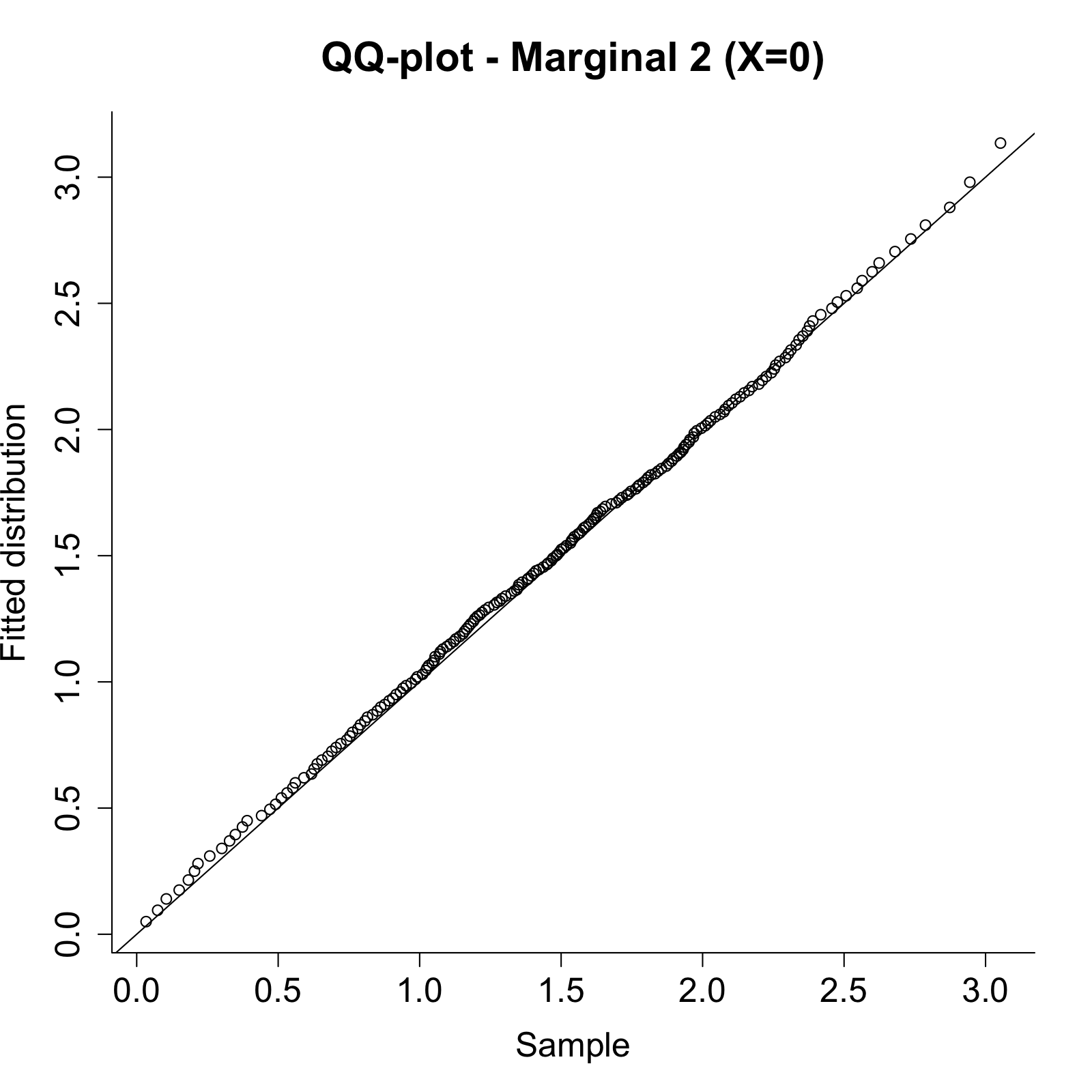}
\includegraphics[width=0.49\textwidth]{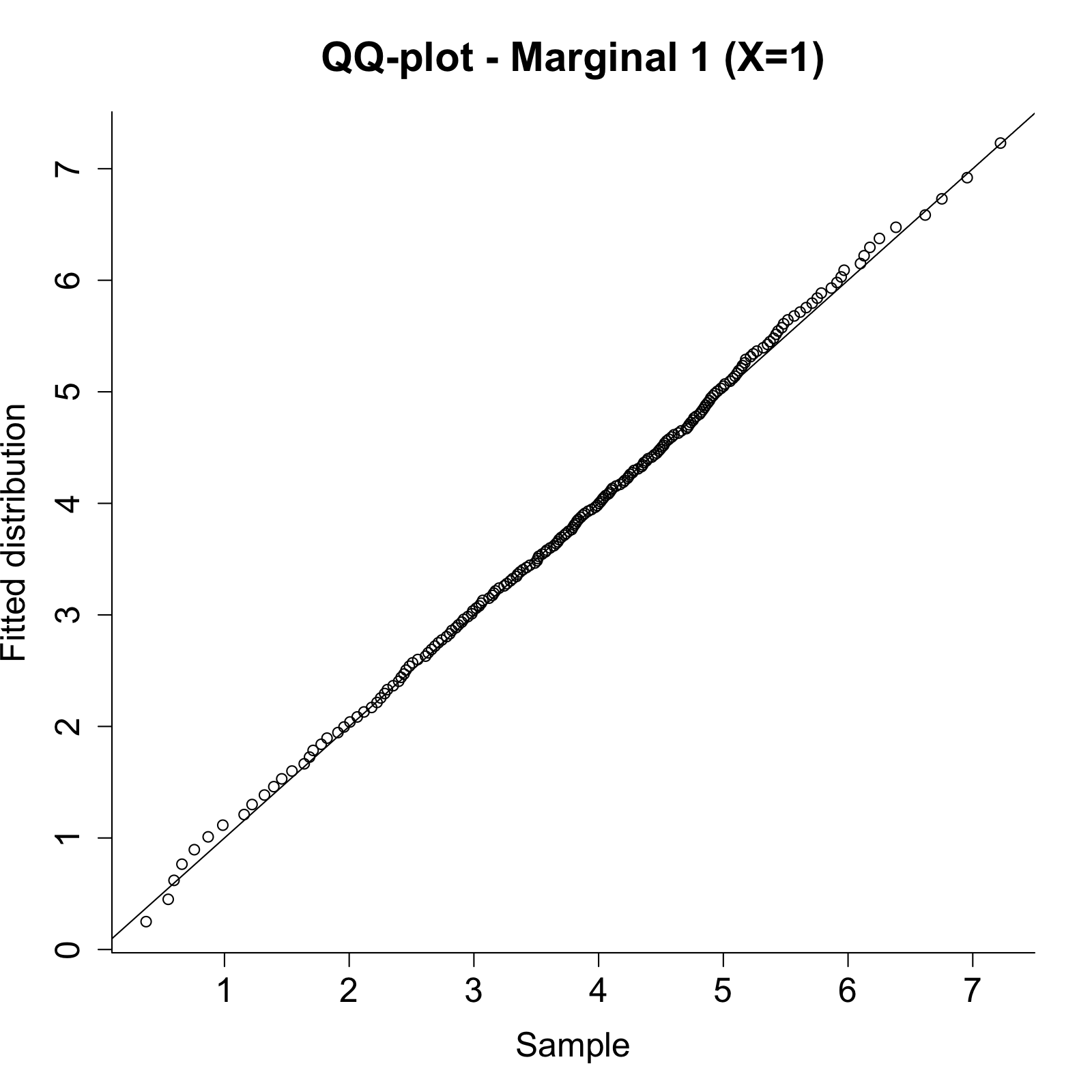}
\includegraphics[width=0.49\textwidth]{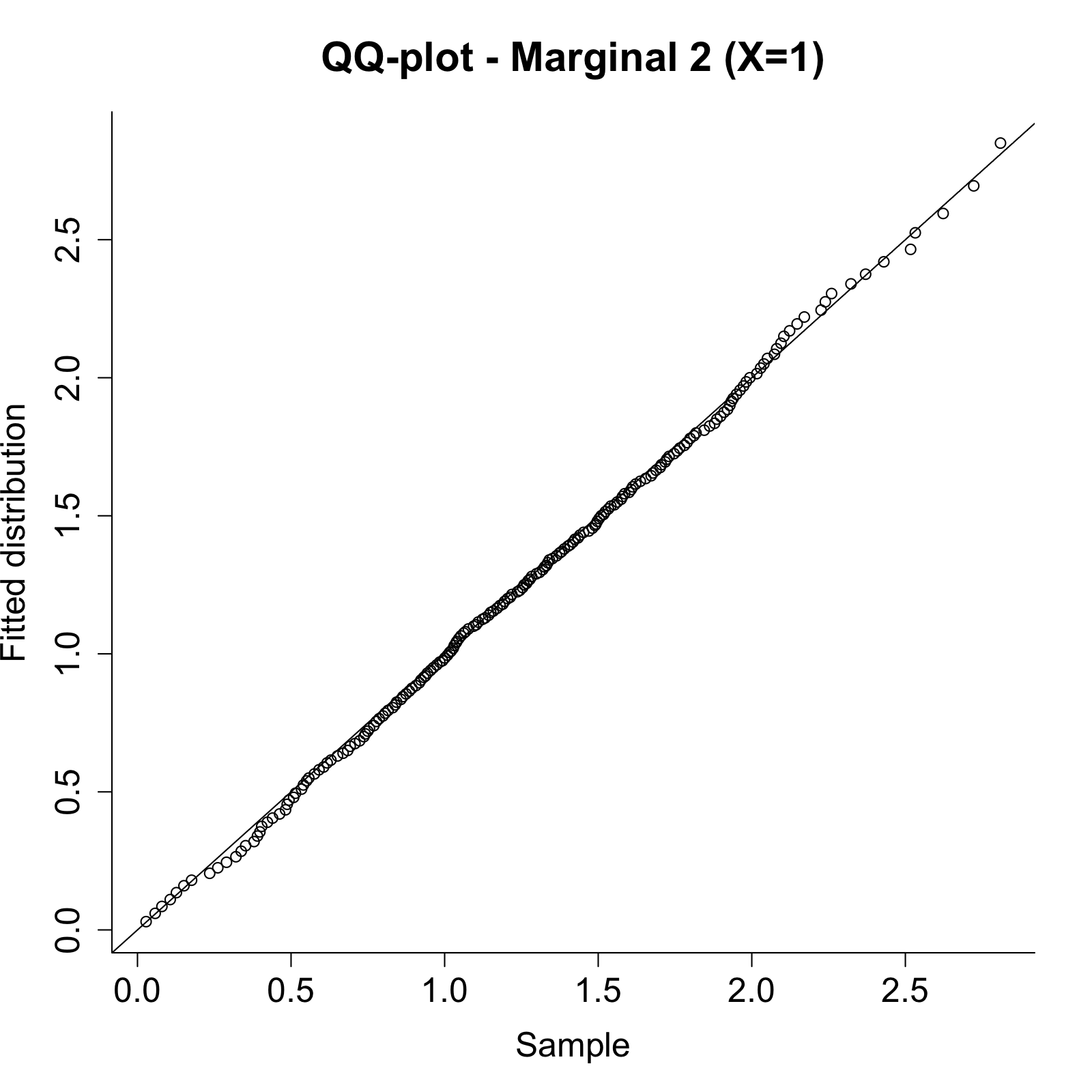}
\caption{QQ-plots of simulated sample from the shared lognormal frailty model versus fitted shared phase-type frailty model.}
\label{fig:ln}
\end{figure}

\begin{figure}[h]
\centering
\includegraphics[width=0.49\textwidth]{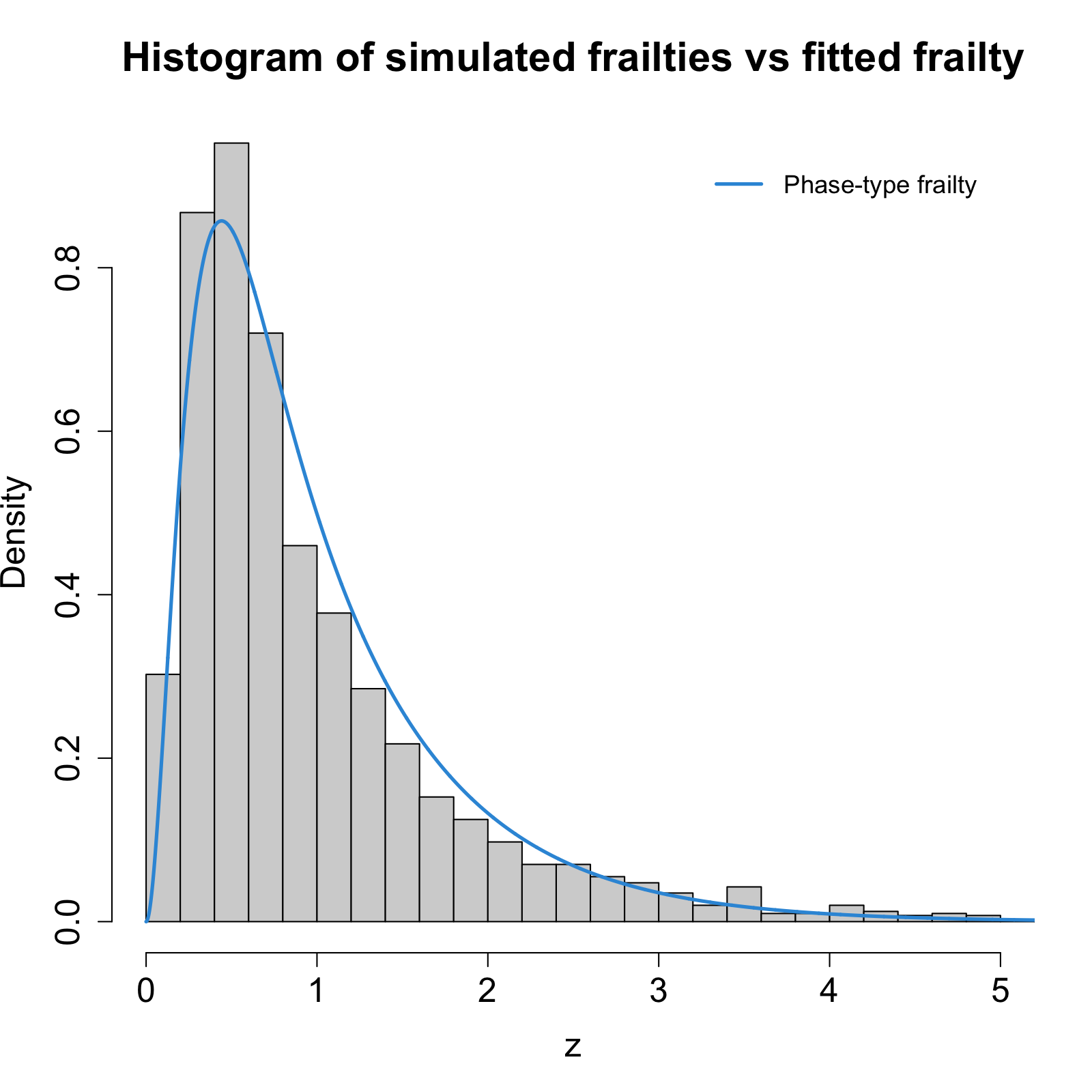}
\caption{Histogram of simulated lognormal frailties versus density of the underlying phase-type frailty (scaled to have mean one).}
\label{fig:ln_under}
\end{figure}

Then, we fit a shared phase-type frailty model of dimension $3$ with Gompertz baseline hazards to the simulated bivariate data. The estimated parameters are
\begin{gather*} 
	\hat{\bfpi}=\left(
	1, \,0, \, 0 \right)\,, \\ 
	\hat{\bfT}=\left( \begin{array}{cccc}
	-1.7179 & 1.7179 & 0.0206 \\
	0.1679 & -1.6930 & 1.5251 \\
	0.5181 & 0.2911 & -1.8967
	\end{array} \right) \,, \\ 
	\hat{b}_1= 0.0034\,, \quad \hat{c}_1= 0.9642\,, \\
	\hat{b}_2= 0.0323\,, \quad \hat{c}_2= 1.9014\,,  \\
	 \hat{\beta}= 0.5485 \,.
\end{gather*}
The quality of the fit is supported by Figure~\ref{fig:ln}, which shows that the algorithm recovers the data structure for both groups.  %{\red 
Additionally, Figure~\ref{fig:ln_under} shows that the density of the underlying phase-type frailty approximates well the histogram of the simulated lognormal frailties. %}

%Alternatively, given a theoretical lognormal (or any other) frailty model, we can approximate such a model by simply fitting a phase-type distribution to the theoretical given distribution. We refer to  \cite{asmussen1996fitting} for details and illustrations, which include the lognormal distribution.% We refer to such reference for further details. 

%{\red
\subsection{Diagnosis of fracture healing} 
\begin{figure}[h]
\centering
\includegraphics[width=0.49\textwidth]{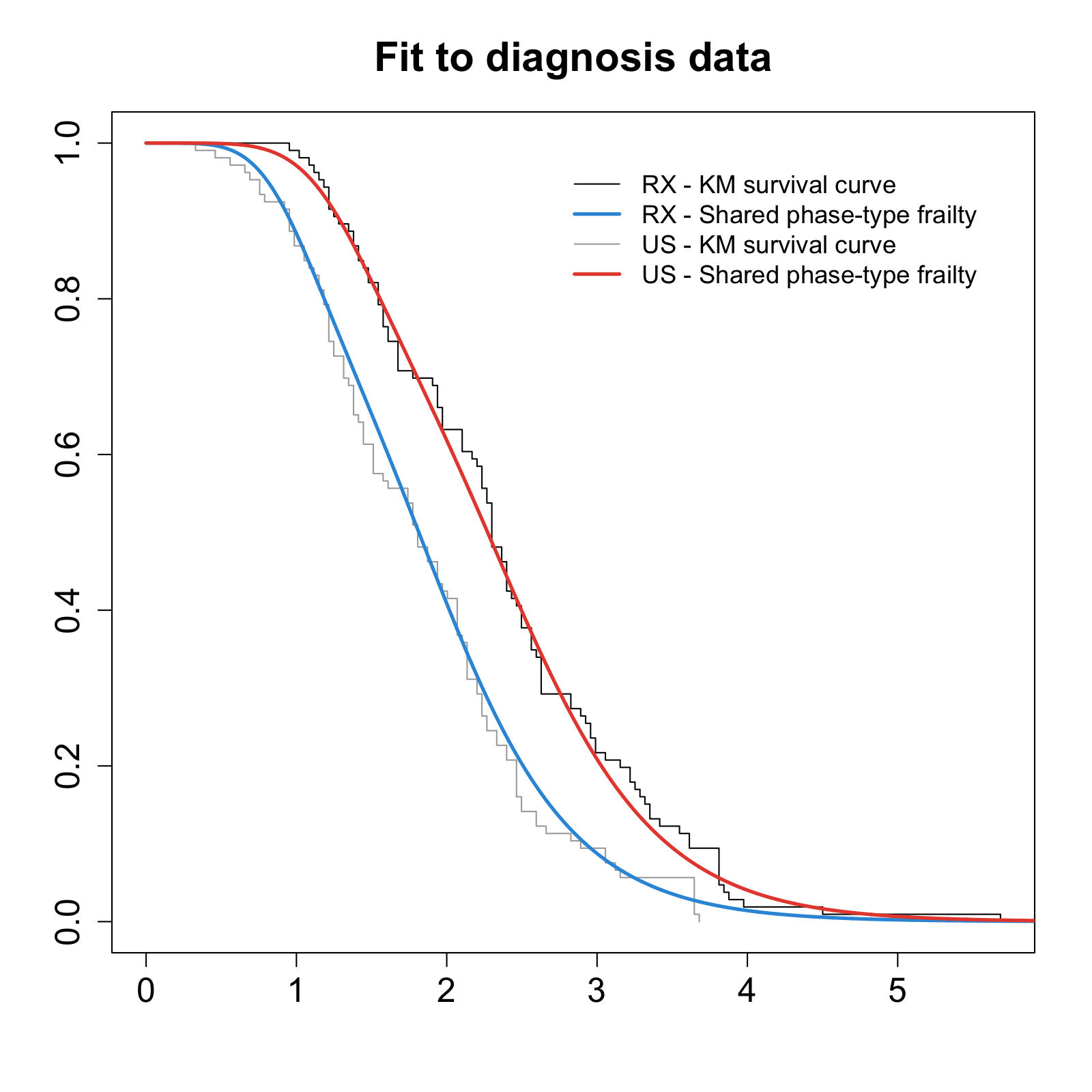}
\includegraphics[width=0.49\textwidth]{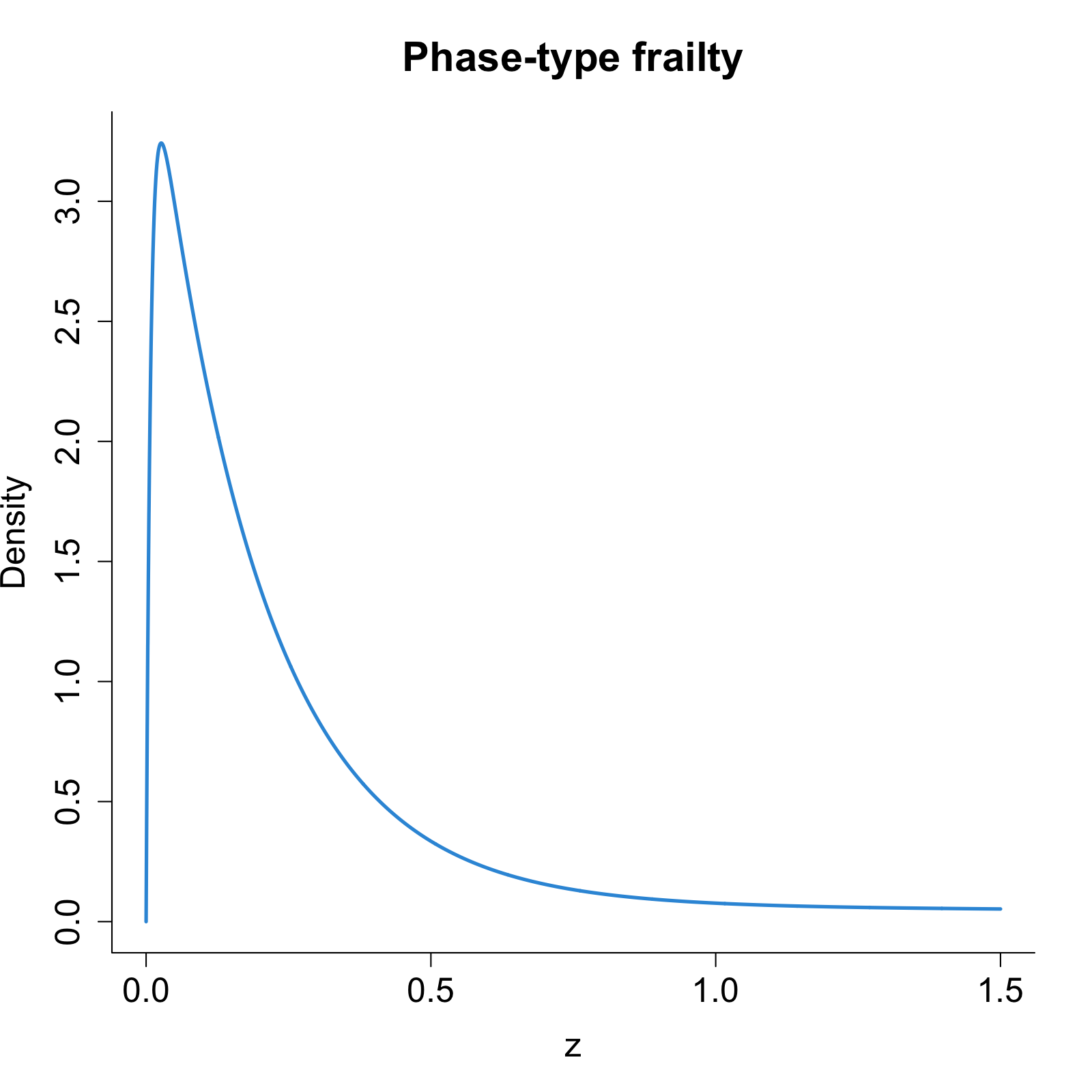}
\caption{Kaplan-Meier survival curves of RX and US versus marginal survival functions from fitted shared phase-type frailty model (left) and density function of the underlying phase-type frailty (right) }
\label{fig:diagnosis}
\end{figure}

We consider the \textit{diagnosis} dataset available in the \texttt{parfm} R-package \citep{munda2012parfm}, which contains information on 106 dogs treated at the veterinary university hospital of Ghent. For each dog, the time to fracture healing was evaluated using two medical imaging techniques: radiography (RX) and ultrasound (US). This dataset was previously analyzed in \cite{duchateau2008frailty} using shared frailty models to jointly model the healing times evaluated under the two techniques, particularly employing Gamma and inverse Gaussian frailties. Here, we fit a shared phase-type frailty model to this data. Specifically, we consider Weibull baseline hazards of the form $\mu_j(y) = \lambda_j \theta_j y^{\theta_j-1}$, $\lambda_j, \theta_j> 0$, $j = 1,2$ and a phase-type frailty of dimension 3 with general Coxian structure of the parameters, obtaining the following estimated parameters
\begin{gather*} 
	\hat{\bfpi}=\left(
	0.2969, \,0.7031, \, 0 \right)\,, \\ 
	\hat{\bfT}=\left( \begin{array}{cccc}
	-0.2346  &  0.2346 &  0 \\
	0 & -124.8573 & 124.8573 \\
	0 &  0 & -5.2942
	\end{array} \right) \,, \\ 
	\hat{\lambda}_1= 0.1131\,, \quad \hat{\theta}_1= 5.0582\,, \\
	\hat{\lambda}_2= 0.0214\,, \quad \hat{\theta}_2= 5.6986\,.  \\
\end{gather*}
Figure~\ref{fig:diagnosis} shows that the fitted model successfully recovers the marginal behavior of the data, with the estimated survival functions aligning closely with the non-parametric Kaplan-Meier survival curves. Furthermore, note that the model leads to a larger value of the log-likelihood ($-219.0$) compared to shared frailty models with  Gamma and inverse Gaussian frailties, which have log-likelihoods of -232.1 and -222.4, respectively. 

%}

\section{Conclusion}\label{sec:conclusion}
We proposed a novel frailty model where phase-type distributions are employed as frailties. We showed that this model exhibits similar properties to the Gamma frailty model with closed-form expressions for different functionals. Moreover, we provided some multivariate extensions, namely the shared frailty and correlated frailty models. We derived EM algorithms for maximum-likelihood estimation for all the proposed models, which were then employed in some numerical illustrations. The results show that phase-type frailty models are interesting modeling tools in survival analysis and insurance applications.  A possible line of work for future research is on generalizing further the correlated phase-type frailty model by employing the class of multivariate phase-type distributions introduced in  \cite{kulkarni1989new}. Under this framework, the correlated Gamma frailty model with integer shape parameters introduced in \cite{yashin1995correlated} results as a particular case, thus the relevance of studying such a more general construction. 
 
%\appendix
%\normalsize
%\subsection*{Implementation details}
 
%\newpage
\appendix
%\normalsize
%
\section{Fitting algorithm for bivariate phase-type distributions}
%the
%{\red 
%
%
%
%
%
\begin{algorithm}[]
\caption{EM algorithm for fitting a bivariate phase-type distribution to a given joint density}\label{alg:bivPH}
\begin{algorithmic}
\normalsize 
\State \textit{\textbf{Input}: Bivariate joint density $g(\cdot)$.}\\
\begin{enumerate}[label=\arabic*.]
\item[0.] Initialize with some ``arbitrary" $( \bfeta, \bfT_{11}, \bfT_{12}, \bfT_{22})$.
	\item Calculate:

For $i=1,\dots,p_{1}, $
 \begin{align*}
	\hat{\eta}_{i}& = \int_{\R^2_+} \dfrac{\eta_{i}\bfe_{i}^{\top} \exp({\bfT_{11} y_{1} })\bfT_{12}\exp({\bfT_{22}y_{2}})(-\bfT_{22}) \bfe}{\bfeta \exp({\bfT_{11}y_{1}})\bfT_{12}\exp({\bfT_{22}y_{2}})(-\bfT_{22}) \bfe} g(\bfy)d\bfy \,.
\end{align*}

For $i,l=1,\dots,p_{1}, \, i \neq l $,
\begin{align*}
\hat{t}_{il} = \dfrac{\mathlarger{ \int_{\R^2_+}
 t_{il}  \dfrac{ { \int_{0}^{y_{1}}}\bfeta \exp({\bfT_{11}u })\bfe_{i} \bfe_{l}^{\top} \exp({\bfT_{11}(y_{1}-u )}) \bfT_{12} \exp({\bfT_{22}y_{2}})(-\bfT_{22}) \bfe du } { \bfeta \exp({\bfT_{11}y_{1}})\bfT_{12}\exp({\bfT_{22}y_{2}})(-\bfT_{22}) \bfe }  g(\bfy)d\bfy }}{\mathlarger{\int_{\R^2_+} \dfrac{ {\int_{0}^{y_{1}}} \bfeta \exp({\bfT_{11}u})\bfe_{i} \bfe_{i}^{\top} \exp({\bfT_{11}(y_{1}-u)}) \bfT_{12} \exp({\bfT_{22}y_{2}}) (-\bfT_{22} ) \bfe  du}{\bfeta \exp({\bfT_{11}y_{1}})\bfT_{12}\exp({\bfT_{22}y_{2}})(-\bfT_{22}) \bfe}  g(\bfy)d\bfy }} \,.
\end{align*}

For $i=1,\dots,p_{1},\, l=p_{1}+1,\dots,p$,
\begin{align*}
\hat{t}_{il} =
   \dfrac{ \mathlarger{ \int_{\R^2_+} t_{il}  \dfrac{\bfeta \exp({\bfT_{11}y_{1} }) \bfe_{i}  \bfe_{l-p_{1}}^{\top}
	 \exp({\bfT_{22}y_{2}}) (-\bfT_{22}) \bfe }{\bfeta \exp({\bfT_{11}y_{1}})\bfT_{12}\exp({\bfT_{22}y_{2}})(-\bfT_{22}) \bfe}  g(\bfy)d\bfy}}{\mathlarger{\int_{\R^2_+} \dfrac{ {\int_{0}^{y_{1}}} \bfeta \exp({\bfT_{11}u})\bfe_{i} \bfe_{i}^{\top} \exp({\bfT_{11}(y_{1}-u)}) \bfT_{12} \exp({\bfT_{22}y_{2}}) (-\bfT_{22} ) \bfe du}{\bfeta \exp({\bfT_{11}y_{1}})\bfT_{12}\exp({\bfT_{22}y_{2}})(-\bfT_{22}) \bfe}  g(\bfy)d\bfy }} \,.
\end{align*}

For $i,l=p_{1}+1,\dots,p, \, i\neq l$,
\begin{align*}
\hat{t}_{il} =
  \dfrac{\mathlarger{\int_{\R^2_+} t_{il}  \dfrac{ {\int_{0}^{y_{2}}} \bfeta \exp({\bfT_{11}y_{1}}) \bfT_{12} \exp({\bfT_{22}u }) \bfe_{i-p_{1}}   \bfe_{l-p_{1}}^{\top} \exp({\bfT_{22}(y_{2}- u)})(-\bfT_{22}) \bfe  du } { \bfeta \exp({\bfT_{11}y_{1}})\bfT_{12} \exp({\bfT_{22}y_{2}})(-\bfT_{22}) \bfe }  g(\bfy)d\bfy}}{\mathlarger{ \int_{\R^2_+} \dfrac{ {\int_{0}^{y_{2}}} \bfeta \exp({\bfT_{11}y_{1}})\bfT_{12}\exp({\bfT_{22}u})\bfe_{i-p_{1}} \bfe_{i-p_{1}}^{\top} \exp({\bfT_{22}(y_{2}-u)})(-\bfT_{22}) \bfe  du }{ \bfeta \exp({\bfT_{11}y_{1}})\bfT_{12} \exp({\bfT_{22}y_{2}})(-\bfT_{22}) \bfe }  g(\bfy)d\bfy }} \,.
\end{align*}

Finally, for $i=p_{1}+1,\dots,p$,
\begin{align*}
	\hat{t}_{i}&  = \dfrac{\mathlarger{ \int_{\R^2_+} t_{i}  \dfrac{\bfeta \exp({\bfT_{11}y_{1} }) \bfT_{12} \exp({\bfT_{22}y_{2} }) \bfe_{i-p_{1}}  }{ \bfeta \exp({\bfT_{11}y_{1}})\bfT_{12}\exp({\bfT_{22}y_{2}})(-\bfT_{22}) \bfe  } g(\bfy)d\bfy }}{ \mathlarger{\int_{\R^2_+} \dfrac{ {\int_{0}^{y_{2}}} \bfeta \exp({\bfT_{11}y_{1}})\bfT_{12}\exp({\bfT_{22}u})\bfe_{i-p_{1}} \bfe_{i-p_{1}}^{\top} \exp({\bfT_{22}(y_{2}-u)})(-\bfT_{22}) \bfe  du}{ \bfeta \exp({\bfT_{11}y_{1}})\bfT_{12}\exp({\bfT_{22}y_{2}})(-\bfT_{22}) \bfe }  g(\bfy)d\bfy }}\,,
\end{align*}
and for $i=1,\dots,p$,
\begin{align*}
	\hat{t}_{ii} = - \sum_{l\neq i} \hat{t}_{il} - \hat{t}_i \,.
\end{align*}

Let $\hat{\bfeta} = (\hat{\eta}_{1}, \dots, \hat{\eta}_{p_1})$, $\hat{\bfT}_{11} = \{\hat{t}_{il}\}_{i,l = 1,\dots,p_1}$, $\hat{\bfT}_{12} = \{\hat{t}_{il}\}_{i = 1,\dots,p_1, l =  p_1 + 1,\dots,p}$, $\hat{\bfT}_{22} = \{\hat{t}_{il}\}_{i,l = p_1 + 1,\dots,p}$, and $\hat{\bft} = (\hat{t}_{p_1 + 1}, \dots, \hat{t}_p)^{\top}$.

\item Assing  $ \bfeta:=\hat{\bfeta}$, $\bfT_{11}:=\hat{\bfT}_{11}$, $\bfT_{12}:=\hat{\bfT}_{12}$, $\bfT_{22}:=\hat{\bfT}_{22}$, and GOTO 1 until a stopping rule is satisfied. 

\end{enumerate}
\State \textit{\textbf{Output}: Fitted parameters $( \bfeta, \bfT_{11}, \bfT_{12}, \bfT_{22})$.}
\end{algorithmic}
\end{algorithm}
%}
%

\newpage
\bibliographystyle{apalike}
\bibliography{references.bib}

\end{document}